\documentclass[11pt]{amsart}
\usepackage[utf8]{inputenc}
\usepackage{amsmath,amsfonts,amssymb}
\usepackage{mathrsfs}
\usepackage{mathtools}
\usepackage{microtype}
\mathtoolsset{showonlyrefs} 
\usepackage[colorlinks=true,linkcolor=blue,citecolor=blue,urlcolor=blue]{hyperref}
\allowdisplaybreaks

\numberwithin{equation}{section}

\theoremstyle{plain}
\newtheorem{theorem}{Theorem}
\newtheorem{lemma}{Lemma}
\newtheorem{proposition}{Proposition}
\newtheorem{corollary}{Corollary}

\theoremstyle{definition}

\newtheorem{remark}{Remark}

\begin{document}

\title[On the Principle of Localization on Quasi-Orbits]{On the Principle of Localization on Quasi-Orbits in Crossed Products}

\author{I.~V.~Baibulov}
\address{Saint Petersburg University, 7/9 Universitetskaya nab., St. Petersburg, 199034 Russia}
\email{i.v.baibulov@gmail.com}

\thanks{This work was performed at the Saint Petersburg Leonhard Euler International Mathematical Institute and supported by the Ministry of Science and Higher Education of the Russian Federation (agreement no.~075-15-2025-343).}


\begin{abstract}
Quasi-regular \(C^*\)-dynamical systems with locally closed quasi-orbits are considered. For the corresponding crossed products, a localization principle on quasi-orbits is proven. This principle reduces the study of the primitive spectrum of the original crossed product to the study of the primitive spectra of crossed products, each corresponding to its own \(C^*\)-dynamical system with a single quasi-orbit. Assuming that the group action on the primitive spectrum of the \(C^*\)-algebra in the crossed product is free, except for some fixed subset, a characterization of the spectrum of the crossed product is provided, which allows for a description of the Jacobson topology on it.
\end{abstract}

\maketitle

\section{Introduction}
Let \((\mathcal{A}, G, \alpha)\) be a \(C^*\)-dynamical system, \(\mathcal{A} \rtimes _\alpha G\) the corresponding crossed product, and \(\mathcal{Q}_{\mathcal{A}}\) the quasi-orbit space of the primitive spectrum \(\operatorname{Prim} \mathcal{A}\). Under certain additional assumptions there exists a continuous surjective map \(\sigma \colon \operatorname{Prim} (\mathcal{A}\rtimes _\alpha G) \to \mathcal{Q}_{\mathcal{A}}\), where each fiber \(\sigma ^{-1}(Q)\) over a quasi-orbit \(Q\) can be identified with the primitive spectrum of an auxiliary \(C^*\)-algebra \(\tilde{\mathcal{J}}_Q \rtimes_{\alpha^Q}G\). Thus, the spectrum \(\operatorname{Prim} (\mathcal{A}\rtimes _\alpha G)\) can be represented as a disjoint union of the spectra \(\operatorname{Prim} (\tilde{\mathcal{J}}_Q \rtimes_{\alpha^Q}G)\). The algebras \(\tilde{\mathcal{J}}_Q \rtimes_{\alpha^Q}G\) are simpler objects compared to \(\mathcal{A}\rtimes _\alpha G\), in particular, the quasi-orbit space \(\mathcal{Q}_{\tilde{\mathcal{J}}_Q}\) consists of a single point. The transition from \(\mathcal{A}\rtimes _\alpha G\) to \(\tilde{\mathcal{J}}_Q \rtimes_{\alpha^Q}G\) is called localization at the quasi-orbit \(Q\), and the theorem describing the spectrum \(\operatorname{Prim} (\mathcal{A}\rtimes _\alpha G)\) as the union of the spectra \(\operatorname{Prim} (\tilde{\mathcal{J}}_Q \rtimes_{\alpha^Q}G)\) is referred to as the localization principle on quasi-orbits. The purpose of this work is to strengthen this theorem and to describe the Jacobson topology on \(\operatorname{Prim} (\mathcal{A}\rtimes _\alpha G)\). This can be achieved under certain restrictions on the action \(\alpha\).

We assume that there exists an invariant ideal \(\mathcal{J}_0 \subset \mathcal{A}\) such that the action \(\alpha\) is free on \(\operatorname{Prim} \mathcal{J}_0 \subset \operatorname{Prim} \mathcal{A}\) and trivial on the remaining part \(\operatorname{Prim} (\mathcal{A}/\mathcal{J}_0)\). Roughly speaking, one can say that if a quasi-orbit \(Q\) lies in \(\mathcal{Q}_{\mathcal{J}_0} \subset \mathcal{Q}_{\mathcal{A}}\), then the fiber \(\sigma ^{-1}(Q)\) consists of a single point and \(\sigma\) restricts to a homeomorphism between \(\operatorname{Prim} (\mathcal{J}_0 \rtimes _\alpha G) \subset \operatorname{Prim} (\mathcal{A}\rtimes _{\alpha } G) \) and \(\mathcal{Q}_{\mathcal{J}_0}\). If \(Q \in \mathcal{Q}_{\mathcal{A}/\mathcal{J}_0}\), then \(Q\) is a singleton, and \(\mathcal{Q}_{\mathcal{A}/\mathcal{J}_0}\) coincides with \(\operatorname{Prim} (\mathcal{A}/\mathcal{J}_0)\). It is assumed that each fiber \(\sigma ^{-1}(Q)\) is homeomorphic to the spectrum of some crossed product \(\mathcal{D} \rtimes _\beta G\), and that there exists a homeomorphism between \(\operatorname{Prim} ((\mathcal{A}/\mathcal{J}_0) \rtimes _{\alpha ^{\mathcal{J}_0} } G)\) and \(\operatorname{Prim} (\mathcal{A}/\mathcal{J}_0 ) \times \operatorname{Prim} (\mathcal{D}\rtimes _\beta G)\). In the final Theorem \ref{th:free-and-fixed-action} and Remark \ref{rem:sim-rel-1}, these two homeomorphisms are glued together into a homeomorphism between \(\operatorname{Prim} (\mathcal{A}\rtimes _\alpha G)\) and the quotient space \((\mathcal{Q}_{\mathcal{A}} \times \operatorname{Prim} (\mathcal{D}\rtimes _\beta G))/_\approx\) with respect to a certain equivalence relation.

A statement close to the mentioned localization principle was discussed in \cite[formula~(4.3)]{Echt.Will2008} under the assumption that the \(C^*\)-algebra \(\mathcal{A}\) in the crossed product is a separable GCR-algebra and the group \(G\) is second countable. In \cite[Corollary~7.27]{Echt2017}, the localization principle was established for smooth \(C^*\)-dynamical systems. In \cite[Theorem~2]{Bayb.Sara2024}, the localization principle was proved under the assumption that the quasi-orbit space of a \(C^*\)-dynamical system \((\mathcal{A},G,\alpha )\) is almost Hausdorff; throughout \cite{Bayb.Sara2024} only discrete groups \(G\) were considered, but this fact was not used in the proof of Theorem 2. In the present paper, the localization principle is proved under the weaker conditions that the dynamical system \((\mathcal{A},G,\alpha )\) is quasi-regular and that each quasi-orbit \(Q \in \mathcal{Q}_{\mathcal{A}}\) is locally closed.

The spectrum \(\operatorname{Prim} (\mathcal{A}\rtimes _\alpha G)\) in the case when the group \(G\) acts freely on \(\operatorname{Prim} \mathcal{A}\) has been studied in \cite{Zell1968, Gree1978, Will1981a, Sier2010, Echt.Laca2013,Echt2017,Bayb.Sara2024}. In \cite{Bayb.Sara2024} the spectrum (together with the Jacobson topology) of a concrete crossed product where the group action on the spectrum is free everywhere except for a certain fixed subset was also described. In the present paper this result from \cite{Bayb.Sara2024} is generalized to a certain class of crossed products, which includes, for example, the crossed products arising in \cite{Bayb2025}.

In Section \ref{sec:preliminaries}, the localization principle on quasi-orbits is proved (Theorem \ref{th:orbit-localisation}). Section \ref{sec:fibers} is devoted to proving the continuity, surjectivity, and openness of the map \(\sigma \colon \operatorname{Prim} (\mathcal{A}\rtimes _\alpha G) \to \mathcal{Q}_{\mathcal{A}}\).
In Section \ref{sec:action-free}, known results on the spectrum of crossed products with free action are restated or slightly generalized; this section is primarily a review. A particular case of a crossed product with trivial action is considered in Section \ref{sec:action-fixed}. The main result of the paper (Theorem \ref{th:free-and-fixed-action}) is proved in Section \ref{sec:action-free-and-fixed}.

The author is grateful to O. V. Sarafanov for posing the problem and for fruitful discussions.

\section{The localization principle on quasi-orbits}
\label{sec:preliminaries}
\subsection{\(C^*\)-dynamical system and crossed product}
Let \(\mathcal{A}\) be a \(C^*\)-algebra, \(G\) a locally compact group,
and \(\alpha\colon t\mapsto \alpha_t\) a fixed continuous homomorphism
of the group \(G\) into the group \(\operatorname{Aut} \mathcal{A}\) of automorphisms of the \(C^*\)-algebra
\(\mathcal{A}\), endowed with the topology of pointwise convergence. The triple
\(\left(\mathcal{A},G,\alpha\right)\) is called a \(C^*\)-dynamical system, and the homomorphism \(\alpha\) is called the action of the group
\(G\) on the \(C^*\)-algebra \(\mathcal{A}\). Throughout, by an algebra we mean a \(C^*\)-algebra, by a dynamical system a \(C^*\)-dynamical system, by a homomorphism of algebras a \(*\)-homomorphism, and by an ideal a closed two-sided ideal.

Denote by \(C_c(G, \mathcal{A})\) the set of continuous functions on \(G\) with
compact support and values in \(\mathcal{A}\). The multiplication and involution
operations in \(C_c(G, \mathcal{A})\) are defined by  
\[
(f * g)(t) = \int_G f(s) \alpha_s(g(s^{-1}t)) \, ds,
\]
\[
f^*(t) = \Delta(t^{-1})\alpha_t(f(t^{-1})^*),
\]
where \(ds\) is a fixed left Haar measure on \(G\), and \(\Delta\) is the modular function of the group \(G\). With these operations the set \(C_c(G, \mathcal{A})\) forms a \(*\)-algebra.

Let \(\pi \colon \mathcal{A} \to B(\mathcal{H})\) be a representation of the algebra \(\mathcal{A}\) on a Hilbert space \(\mathcal{H}\), and let \(U \colon G \to \mathcal{U}(\mathcal{H}) \colon t \mapsto U_t\) be a unitary representation of the group \(G\) on the same space \(\mathcal{H}\) (that is, \(U\) is a continuous homomorphism from the group \(G\) into the group \(\mathcal{U}(\mathcal{H})\) of unitary operators on \(\mathcal{H}\), endowed with the topology of pointwise convergence). If for all \(t \in G\) and \(a \in \mathcal{A}\) the relation  
\[
\pi(\alpha_t(a)) = U_t \pi(a) U_t^*
\]  
holds, then the pair \((\pi, U)\) is called a covariant representation of the dynamical system \((\mathcal{A}, G, \alpha)\). A covariant representation \((\pi ,U)\) is called non-degenerate if the representation \(\pi\) is non-degenerate. Every covariant representation \((\pi, U)\) of the system \((\mathcal{A}, G, \alpha)\) gives rise to a representation of the \(*\)-algebra \(C_c(G, \mathcal{A})\):
\begin{equation}\label{eq:integr_rep}
(\pi \rtimes U)(f) = \int_G \pi(f(t)) U_t \, dt, \quad f \in C_c(G, \mathcal{A}).
\end{equation}
The crossed product of \(\mathcal{A}\) and \(G\) is defined as the completion of \(C_c(G, \mathcal{A})\) with respect to the \(C^*\)-norm 
\[
\|f\| = \sup \{ \| \pi \rtimes U (f) \| : (\pi, U) \text{ is a covariant representation of } (\mathcal{A}, G, \alpha) \},
\]  
denoted by \(\mathcal{A} \rtimes_\alpha G\) (see \cite[Chapter~2]{Will2007} for details). Note that in the literature the algebra \(\mathcal{A} \rtimes_\alpha G\) is also called the full crossed product. In this paper we omit the adjective ``full'' for brevity.

The continuous extension of the representation \(\pi\rtimes U\) from \eqref{eq:integr_rep} of the \(*\)-algebra \(C_c(G,\mathcal{A})\) to a representation of the \(C^*\)-algebra \(\mathcal{A}\rtimes_{\alpha}G\) will be denoted by the same symbol \(\pi\rtimes U\).  
By \cite[Proposition~2.40]{Will2007}, the correspondence \((\pi ,U) \mapsto  \pi \rtimes U\) between non-degenerate covariant representations of \((\mathcal{A}, G, \alpha)\) and non-degenerate representations of the algebra \(\mathcal{A}\rtimes_{\alpha}G\) is bijective, and a representation \(\pi\rtimes U\) of the algebra \(\mathcal{A}\rtimes_\alpha G\) is irreducible if and only if the covariant representation \((\pi,U)\) is irreducible (that is, any closed subspace that reduces both \(\pi\) and \(U\) is either \(\{0\}\) or \(\mathcal{H}\)). 

An ideal \(\mathcal{I}\) in the algebra \(\mathcal{A}\) is called invariant if \(\alpha _t(\mathcal{I}) = \mathcal{I}\) for all \(t \in G\). Each invariant ideal \(\mathcal{I}\) gives rise to an ideal \(\mathcal{I}\rtimes_{\alpha}G\) in \(\mathcal{A} \rtimes _\alpha G\), which coincides with the completion of the set \(C_c(G,\mathcal{I}) \subset C_c(G,\mathcal{A})\) with respect to the norm in \(\mathcal{A}\rtimes_{\alpha}G\) (see \cite[Lemma~3.17]{Will2007}); the map \(C_c(G,\mathcal{A}) \to C_c(G,\mathcal{A}/\mathcal{I})\colon f \mapsto \rho _{\mathcal{I}}(f(\cdot ))\), where \(\rho _{\mathcal{I}} \colon \mathcal{A} \to  \mathcal{A}/\mathcal{I}\) is the quotient map, induces an isomorphism
\begin{equation}\label{eq:full-cr-prod-exact}
	(\mathcal{A} \rtimes _\alpha G)/(\mathcal{I} \rtimes _\alpha G) \cong (\mathcal{A}/\mathcal{I}) \rtimes _{\alpha^{\mathcal{I}} } G,
\end{equation}
where \(\alpha^{\mathcal{I}}\) is the action on the quotient \(\mathcal{A}/\mathcal{I}\) induced by the action \(\alpha\) (see \cite[Proposition~3.19]{Will2007}). In what follows, we will not distinguish between the algebras \((\mathcal{A} \rtimes _\alpha G)/(\mathcal{I} \rtimes _\alpha G)\) and \((\mathcal{A}/\mathcal{I}) \rtimes _{\alpha^{\mathcal{I}} } G\), having in mind the above isomorphism.

\subsection{The spectrum of the crossed product}
\label{subsec:crossed-prod-spec}
We recall some notions necessary for studying the spectrum of a crossed product. The action \(\alpha\) of the group \(G\) induces a jointly continuous action on the primitive spectrum \(\operatorname{Prim} \mathcal{A}\), endowed with the Jacobson topology, by the rule \(P \mapsto \alpha_t(P)\) (see \cite[Lemma~1.3]{Glim1962} or \cite[Lemma~7.1]{Raeb.Will1998}). The action thus obtained will also be denoted by \(\alpha\). Primitive ideals \(P_1\) and \(P_2\) are said to belong to the same quasi-orbit if
\begin{equation}\label{eq:quasi-orbit-sim}
	\overline{\{\alpha_t(P_1): t \in G\}} = \overline{\{\alpha_t(P_2): t \in G\}};
\end{equation}
this is an equivalence relation, denoted by \(\sim_q\); the corresponding quotient space
\(\mathcal{Q}_{\mathcal{A}}=\operatorname{Prim}\mathcal{A} /{\sim_q}\)
is called the quasi-orbit space of the dynamical system
\((\mathcal{A}, G, \alpha)\). Note that by \cite[Lemma~6.12]{Will2007} the quotient map \(\kappa\colon \operatorname{Prim} \mathcal{A} \to \mathcal{Q}_{\mathcal{A}} \) is open, hence the quasi-orbit space is a continuous open image of a locally compact space and therefore is itself locally compact.

For every quasi-orbit \(Q\in \mathcal{Q}_{\mathcal{A}}\), we define the associated invariant ideal
\(\mathcal{I}_Q =\bigcap_{P \in \kappa^{-1}(Q)} P \). The action \(\alpha\) induces an action \(\alpha^Q\) on the quotient algebra \(\mathcal{A}_Q =
\mathcal{A}/\mathcal{I}_Q\). Since \(\mathcal{A}_Q
\rtimes_{\alpha^Q}G\) is a quotient of \(\mathcal{A}
\rtimes_{\alpha} G\), we identify the spectrum \(\left(\mathcal{A}_Q
\rtimes_{\alpha^Q}G\right)\widehat{}\,\) with the corresponding closed subset of \(\left(\mathcal{A}
\rtimes_{\alpha} G\right)\widehat{}\,\),
so that
\[\left(\mathcal{A} \rtimes_{\alpha} G\right)\widehat{}  \supset \bigcup_{Q \in
\mathcal{Q}_{\mathcal{A}}} \left(\mathcal{A}_Q
\rtimes_{\alpha^Q}G\right)\widehat{}\,.\] 
A dynamical system
\((\mathcal{A},G,\alpha )\) is called quasi-regular if for every irreducible
representation \(L= \pi \rtimes U\) of \(\mathcal{A}\rtimes_\alpha G\)
there exists a quasi-orbit \(Q_\pi \in \mathcal{Q}_{\mathcal{A}}\) such that
\begin{equation}\label{Eq:lives-on-quasi-orbit}
	\ker \pi = \mathcal{I}_{Q_\pi}.
\end{equation}
In other words, the representation \(\pi\) of \(\mathcal{A}\) induces a representation \(\tilde{\pi }\) of the quotient \(\mathcal{A}_{Q_\pi }\), and hence \(\tilde{\pi }\rtimes U\) is an irreducible representation of \(\mathcal{A}_{Q_\pi } \rtimes _{\alpha^{Q_\pi }} G\). Therefore, we have the following statement.
\begin{lemma}\label{lm:localization-A-Q}
	If the \(C^*\)-dynamical system
	\((\mathcal{A},G,\alpha )\) is quasi-regular, then the following relations hold:
\begin{align}
	\label{eq:cr-prod_spectrum1:1}
	\left(\mathcal{A} \rtimes_{\alpha} G\right)\widehat{} & = \bigcup_{Q \in \mathcal{Q}_{\mathcal{A}}} \left(\mathcal{A}_Q \rtimes_{\alpha^Q}G\right)\widehat{}, \\
	\label{eq:cr-prod_spectrum1:2}
	\operatorname{Prim} \left(\mathcal{A} \rtimes_{\alpha} G\right) &= \bigcup_{Q \in \mathcal{Q}_{\mathcal{A}}} \operatorname{Prim} (\mathcal{A}_Q \rtimes_{\alpha^Q}G).
\end{align}
\end{lemma}

A topological space is called almost Hausdorff if it is locally compact and every nonempty locally compact subspace contains a nonempty relatively open Hausdorff subspace.
\begin{lemma}[see {\cite[Corollary~19]{Gree1978}}]
\label{lm:quasi-regularity}
	If the quasi-orbit space is second countable or almost Hausdorff, then the \(C^*\)-dynamical system \((\mathcal{A},G,\alpha )\) is quasi-regular.
\end{lemma}
\begin{remark}\label{rem:separable-a-2nd_count-Q}
If the \(C^*\)-algebra
\(\mathcal{A}\) is separable, then the primitive spectrum \(\operatorname{Prim} \mathcal{A}\) is second countable (see \cite[3.3.4]{Diks1974}). Hence, the quasi-orbit space \(\mathcal{Q}_{\mathcal{A}}\) of the \(C^*\)-dynamical system \((\mathcal{A},G,\alpha )\) is second countable.
\end{remark}

The union on the right-hand side of equalities \eqref{eq:cr-prod_spectrum1:1}--\eqref{eq:cr-prod_spectrum1:2} is, generally speaking, not disjoint,
but it can be ``thinned out''. To this end, consider the subset \(R_Q\)
of \(\operatorname{Prim} \left(\mathcal{A}_Q
\rtimes_{\alpha^Q}G\right)\), consisting of those primitive ideals \(\ker (\pi
\rtimes U)\) such that \(\ker \pi = \mathcal{I}_Q \), that is, \(\pi \) is a faithful
representation of the algebra \(\mathcal{A}_Q\). Note that from the inclusion \(\ker (\pi \rtimes U) \in R_{Q_1}\cap R_{Q_2}\) it follows that \(\mathcal{I}_{Q_1} = \ker \pi = \mathcal{I}_{Q_2}\), hence \(Q_1 = Q_2\). Therefore, the union \(\bigcup_{Q \in \mathcal{Q}_{\mathcal{A}}} R_Q
\) is disjoint.

Suppose that a quasi-orbit \(Q \in \mathcal{Q}_{\mathcal{A}}\) is locally closed, that is, the set \(\{ Q \}\) is relatively open in its closure. Then, the set \(R_Q\) is the spectrum of a certain crossed product. To describe such a crossed product, put \(\mathcal{J}_Q= \bigcap_{P\in \overline{\kappa^{-1}(Q)}\setminus\kappa^{-1}(Q)}P\); in the case \(\overline{\kappa^{-1}(Q)}\setminus\kappa^{-1}(Q) = \varnothing\) we put \(\mathcal{J}_Q = \mathcal{A}\). Since the map \(\kappa \) is open, the set \(\overline{\kappa^{-1}(Q)}\setminus\kappa^{-1}(Q)\) is closed in \(\operatorname{Prim} \mathcal{A}\). By the natural bijection between closed ideals of the algebra \(\mathcal{A}\) and open subsets of \(\operatorname{Prim} \mathcal{A}\), the ideal \(\mathcal{I}_Q\) corresponds to the set \(\operatorname{Prim} \mathcal{A} \setminus \overline{\kappa^{-1}(Q)}\), and the ideal \(\mathcal{J}_Q\) corresponds to the set \(\operatorname{Prim} \mathcal{A} \setminus (\overline{\kappa^{-1}(Q)}\setminus\kappa^{-1}(Q))\). Since \(\overline{\kappa^{-1}(Q)} \neq \overline{\kappa^{-1}(Q)}\setminus\kappa^{-1}(Q)\), we have \(\mathcal{I}_Q\neq \mathcal{J}_Q\), and the quotient algebra \(\tilde{\mathcal{J}}_Q= \mathcal{J}_Q/\mathcal{I}_Q\) is nontrivial. 

The algebra \(\tilde{\mathcal{J}}_Q\) is canonically identified with the image of the ideal \(\mathcal{J}_Q\) under the quotient map \(\mathcal{A}\to \mathcal{A}_Q = \mathcal{A}/\mathcal{I}_Q\), and its primitive spectrum \(\operatorname{Prim} \tilde{\mathcal{J}}_Q\) is identified with the subset \(\kappa ^{-1}(Q) \subset \operatorname{Prim} \mathcal{A}\). By \cite[Lemma~3.17]{Will2007} the spectrum \((\tilde{\mathcal{J}}_Q \rtimes_{\alpha^Q}G)\thinspace\widehat{} \,\) is canonically identified with a part of the spectrum \(\left(\mathcal{A}_Q \rtimes_{\alpha^Q}G\right)\widehat{}\), which we, in turn, regard as a part of \(\left(\mathcal{A} \rtimes_{\alpha} G\right)\widehat{}\). 
\begin{lemma}\label{lm:locally-closed-R-Q}
	Let \((\mathcal{A}, G, \alpha)\) be an arbitrary \(C^*\)-dynamical system.
	If a quasi-orbit \(Q \in \mathcal{Q}_{\mathcal{A}}\) is locally closed, then \(R_Q = \operatorname{Prim} (\tilde{\mathcal{J}}_Q \rtimes_{\alpha^Q}G)\).
\end{lemma}
\begin{proof}
	The statement of the lemma was essentially established in the proof of \cite[Theorem~2]{Bayb.Sara2024}. For completeness, we provide a full proof.

	Let \(L=\pi\rtimes U\) be an irreducible representation of \(\mathcal{A}_Q\rtimes_{\alpha^Q}G\). If \(\ker L\notin\operatorname{Prim}(\tilde{\mathcal{J}}_Q\rtimes_{\alpha^Q}G)\), then \(\ker L\supset\tilde{\mathcal{J}}_Q\rtimes_{\alpha^Q}G\), which implies \(\tilde{\mathcal{J}}_Q\subset\ker\pi\). Since \(\tilde{\mathcal{J}}_Q\) is nontrivial, \(\pi\) is not faithful on \(\mathcal{A}_Q\), and hence \(\ker L\notin R_Q\). Thus, \(R_Q\subset\operatorname{Prim}(\tilde{\mathcal{J}}_Q\rtimes_{\alpha^Q}G)\).

	Conversely, let \(L=\pi\rtimes U\) be an irreducible representation of \(\mathcal{A}\rtimes_\alpha G\) with \(\ker L\in\operatorname{Prim}(\tilde{\mathcal{J}}_Q\rtimes_{\alpha^Q}G)\). We show that \(\ker\pi=\mathcal{I}_Q\). Let \(K\subset\operatorname{Prim}\mathcal{A}\) be the set of all primitive ideals containing \(\ker\pi\). Since \(\pi\) is nontrivial, \(K\) is nonempty. Moreover, \(K\) is a closed set invariant under \(\alpha\) and \(\ker\pi=\bigcap_{P\in K}P\). Since \(\pi\) corresponds to a representation of \(\tilde{\mathcal{J}}_Q=\mathcal{J}_Q/\mathcal{I}_Q\), we have \(\mathcal{I}_Q\subset\ker\pi\). Thus, \(K\subset\overline{\kappa^{-1}(Q)}\). Suppose \(\ker\pi\neq\mathcal{I}_Q\), so that \(K\neq\overline{\kappa^{-1}(Q)}\). Then, the sets \(\kappa^{-1}(Q)\) and \(K\) are disjoint; otherwise, by the invariance of \(K\), \(\overline{\kappa^{-1}(Q)}\subset K\), implying that \(\overline{\kappa^{-1}(Q)}=K\). Therefore, \(K\subset\overline{\kappa^{-1}(Q)}\setminus\kappa^{-1}(Q)\), that is, \(\ker\pi\supset\mathcal{J}_Q\). This means that the corresponding representation of \(\tilde{\mathcal{J}}_Q\) is trivial. The contradiction obtained proves that \(\ker\pi=\mathcal{I}_Q\), i.e. \(\ker L\in R_Q\).
\end{proof}
\begin{remark}\label{rem:localization-R-Q-A-Q}
	If a quasi-orbit \(Q \in \mathcal{Q}_{\mathcal{A}}\) is closed, then \(R_Q = \operatorname{Prim} \left(\mathcal{A}_Q \rtimes_{\alpha^Q}G\right)\).
\end{remark}
\begin{theorem}[Localization principle on quasi-orbits]\label{th:orbit-localisation}
	Assume that every quasi-orbit \(Q \in \mathcal{Q}_{\mathcal{A}}\) of a quasi-regular \(C^*\)-dynamical system \((\mathcal{A}, G, \alpha)\) is locally closed. Then, we have
	\begin{align}
		\label{Eq:cr-prod_spectrum2:1}
		\left(\mathcal{A} \rtimes_{\alpha} G\right)\widehat{} & = \bigcup_{Q \in \mathcal{Q}_{\mathcal{A}}} (\tilde{\mathcal{J}}_Q \rtimes_{\alpha^Q}G)\thinspace\widehat{} \ ,\\
		\label{Eq:cr-prod_spectrum2:2}
		\operatorname{Prim} \left(\mathcal{A} \rtimes_{\alpha} G\right) &= \bigcup_{Q \in \mathcal{Q}_{\mathcal{A}}} \operatorname{Prim} (\tilde{\mathcal{J}}_Q \rtimes_{\alpha^Q}G),
	\end{align}
	and the unions are disjoint.
\end{theorem}
\begin{proof}
	It is enough to combine Lemmas~\ref{lm:localization-A-Q} and \ref{lm:locally-closed-R-Q}.
\end{proof}

\section{Fiber structure of the spectrum of the crossed product over the quasi-orbit space}
\label{sec:fibers} 

Let \((\mathcal{A},G,\alpha )\) be a quasi-regular dynamical system. Then, \(\operatorname{Prim} (\mathcal{A}\rtimes _\alpha G)\) is the union of the disjoint sets \(R_Q\) (introduced after Remark~\ref{rem:separable-a-2nd_count-Q}). Consider the map \(\sigma \colon \operatorname{Prim} (\mathcal{A}\rtimes _\alpha G) \to \mathcal{Q}_{\mathcal{A}}\) such that \(\sigma ^{-1} (Q) = R_Q\). In this section we study the properties of the map \(\sigma\); in the following sections these results will allow us, under some additional assumptions, to strengthen Theorem~\ref{th:orbit-localisation}, in particular, to describe the topology on \(\operatorname{Prim} (\mathcal{A}\rtimes _\alpha G)\).

\begin{proposition}\label{prop:prim_is_fibered}
	If \((\mathcal{A},G,\alpha )\) is a quasi-regular \(C^*\)-dynamical system, then
	the map \( \sigma \colon \operatorname{Prim} (\mathcal{A}\rtimes _{\alpha }G) \to \mathcal{Q}_{\mathcal{A}}\) is continuous, and its range is dense in \(\mathcal{Q}_{\mathcal{A}} \).
\end{proposition}
\begin{proof}
	If the closure of the image \(\sigma(\operatorname{Prim} (\mathcal{A}\rtimes _{\alpha }G)) \) does not coincide with \( \mathcal{Q}_{\mathcal{A}}\), then the set \(\kappa ^{-1}(\overline{\sigma(\operatorname{Prim} (\mathcal{A}\rtimes _{\alpha }G))}) \) is closed and does not coincide with \(\operatorname{Prim} \mathcal{A}\). Therefore, the invariant ideal 
	\[
	\mathcal{I} = \bigcap_{Q \in \overline{\sigma(\operatorname{Prim} (\mathcal{A}\rtimes _{\alpha }G))} } \mathcal{I}_Q = \bigcap_{P \in \kappa ^{-1}(\overline{\sigma(\operatorname{Prim} (\mathcal{A}\rtimes _{\alpha }G))}) } P 
	\]
	is nontrivial. Hence, the primitive spectrum \(\operatorname{Prim} (\mathcal{I}\rtimes _\alpha G)\) is nonempty and does not intersect the preimage \(\sigma ^{-1}(\overline{\sigma(\operatorname{Prim} (\mathcal{A}\rtimes _{\alpha }G))}) \subset  \operatorname{Prim}( (\mathcal{A}/\mathcal{I}) \rtimes _\alpha G)\). This contradicts the equality \(\sigma ^{-1}(\overline{\sigma(\operatorname{Prim} (\mathcal{A}\rtimes _{\alpha }G))}) = \operatorname{Prim}(\mathcal{A} \rtimes _\alpha G)\). Therefore, the image of the map \(\sigma \) is dense in \(\mathcal{Q}_{\mathcal{A}} \). It remains to prove the continuity of \(\sigma \).

	Recall that for an arbitrary \(C^*\)-algebra \(\mathcal{A}\), the Fell topology on the space of all its closed two-sided ideals \(\operatorname{I}(\mathcal{A})\) is the topology generated by the subbase \(\{ U_{\mathcal{J}} \colon \mathcal{J} \in \operatorname{I} (\mathcal{A})\} \), where
	\begin{equation} \label{Eq:fell_prebase-system}
		U_{\mathcal{J}}=\{\mathcal{I} \in \operatorname{I}(\mathcal{A}): \mathcal{I} \not \supset \mathcal{J} \}.
	\end{equation}
	Consider the map 
	\begin{equation}\label{eq:def-res}
		\operatorname{Res} \colon
	\operatorname{I}(\mathcal{A}\rtimes_{\alpha} G) \to  \operatorname{I}(\mathcal{A}) \colon \operatorname{Res}(\ker L) = \ker \pi,
	\end{equation}
	where \(L = \pi \rtimes U\) is a representation of the algebra \(\mathcal{A}\rtimes _\alpha G\). By \cite[Proposition~9]{Gree1978}, the map \(\operatorname{Res}\) is continuous.
	The Jacobson topology on the set of primitive ideals of \(\mathcal{A}\rtimes_{\alpha} G\) coincides with the topology induced by the Fell topology, that is, \(\operatorname{Prim} (\mathcal{A}\rtimes _\alpha G)\) is a subspace of \(\operatorname{I}(\mathcal{A}\rtimes_{\alpha} G)\). Hence, the restriction of \(\operatorname{Res}\) to \(\operatorname{Prim} (\mathcal{A}\rtimes _\alpha G)\) is also continuous.

	Further, consider the set \(\operatorname{Prim} ^G(\mathcal{A})= \{\mathcal{I}_{Q} \colon Q \in \mathcal{Q}_{\mathcal{A}} \}\) with the topology induced by the Fell topology on \(\operatorname{I}(\mathcal{A}) \).
	By \cite[Lemma on p.~221]{Gree1978}, the map \( \mathcal{Q}_{\mathcal{A}} \to \operatorname{Prim} ^G(\mathcal{A}) \colon  Q \mapsto \mathcal{I}_Q\) is a homeomorphism. Denote the inverse map \(\mathcal{I}_Q \mapsto Q\) by \(q\).

	By the quasi-regularity of the dynamical system \((\mathcal{A},G,\alpha )\), the image \(\operatorname{Res}(\operatorname{Prim} (\mathcal{A}\rtimes _\alpha G))\) is contained in \(\operatorname{Prim} ^G(\mathcal{A})\). Therefore, the composition \(q \circ \operatorname{Res} \colon \operatorname{Prim} (\mathcal{A}\rtimes _\alpha G) \to \mathcal{Q}_{\mathcal{A}}\) is well-defined and continuous. At the same time, its value on the primitive ideals from the set \( R_Q\) is equal to \(Q\). Hence, the map \(\sigma \) coincides with \(q \circ \operatorname{Res}\), and therefore is continuous.

	The proposition is proved.
\end{proof}

\begin{lemma}\label{lm:sigma-surjective}
	Let every quasi-orbit \(Q \in \mathcal{Q}_{\mathcal{A}}\) of a quasi-regular \(C^*\)-dynamical system \((\mathcal{A},G,\alpha )\) be locally closed. Then, the map \( \sigma \colon \operatorname{Prim} (\mathcal{A}\rtimes _{\alpha }G) \to \mathcal{Q}_{\mathcal{A}}\) is surjective. 
\end{lemma}
\begin{proof}
	By Lemma \ref{lm:locally-closed-R-Q}, for every \(Q \in \mathcal{Q}_{\mathcal{A}}\) the set \(\sigma ^{-1}(Q)= R_Q = \operatorname{Prim} (\tilde{\mathcal{J}}_Q \rtimes_{\alpha^Q}G)\) is nonempty.
\end{proof}

Now we study conditions under which the map \(\sigma \colon \operatorname{Prim} (\mathcal{A}\rtimes _{\alpha }G) \to \mathcal{Q}_{\mathcal{A}}\) is open. For this we will need properties of the map \(\operatorname{Ind} \colon \operatorname{I} (\mathcal{A} )\to \operatorname{I}(\mathcal{A}\rtimes _\alpha G) \) constructed by Green in \cite{Gree1978}; see also \cite[Chapter~5]{Will2007} or \cite[Section~7]{Echt2017}. To be precise, these works consider a much more general construction, whose explicit description is lengthy. However, the special case we need admits a simple description. Let \(\mathcal{I} = \ker \pi \in \operatorname{I}( \mathcal{A} )\), where \(\pi\) is a nondegenerate representation of \(\mathcal{A}\) on a Hilbert space \(\mathcal{H}\). Consider the representation \(\tilde{\pi}\) of \(\mathcal{A}\) on the Hilbert space \(L_2(G, \mathcal{H})\) acting by \((\tilde{\pi }(a)\xi )(t)=\pi (\alpha _{t^{-1}}(a))(\xi (t))\), where \(a \in  \mathcal{A}\), \(t \in G\) and \(\xi \in L_2(G, \mathcal{H})\). Denote by \(\lambda\) the unitary representation of \(G\) on \(L_2(G, \mathcal{H})\) given by \((\lambda (s)\xi )(t) = \xi (s^{-1} t)\) for all \(s,t \in G\) and \(\xi \in L_2(G, \mathcal{H})\). 
Then, \(\operatorname{Ind} \mathcal{I} = \ker \tilde{\pi} \rtimes \lambda\). (See \cite[Remark~7.5(3)]{Echt2017}.)  
From \cite[Proposition~9]{Gree1978} it follows that the map \(\operatorname{Ind} \colon \operatorname{I}(\mathcal{A}) \to \operatorname{I}(\mathcal{A}\rtimes _\alpha G )\) preserves intersections of ideals. From \cite[Proposition~13]{Gree1978} it follows that if \(G\) is amenable, then for invariant ideals \(\mathcal{I} \subset \mathcal{A}\) one has \(\operatorname{Ind} \mathcal{I} = \mathcal{I}\rtimes _\alpha G\). 

\begin{proposition}\label{prop:prim_is_fibered:opennes}
	Let \((\mathcal{A},G,\alpha )\) be a quasi-regular \(C^*\)-dynamical system, and let \(G\) be amenable. Then, the following conditions are equivalent.
	\begin{enumerate}
		\item The map \(\sigma \colon \operatorname{Prim} (\mathcal{A}\rtimes _{\alpha }G) \to \mathcal{Q}_{\mathcal{A}}\) is relatively open.
		\item For every \(Q \in \sigma (\operatorname{Prim}(\mathcal{A}\rtimes _{\alpha }G))\) one has \(\overline{\sigma ^{-1} (Q)}  = \sigma ^{-1}(\overline{\{ Q \}})\).
		\item For every \(Q \in \sigma (\operatorname{Prim}(\mathcal{A}\rtimes _{\alpha }G))\) the set \(R_Q\) is dense in \( \operatorname{Prim} (\mathcal{A}_{Q}\rtimes _{\alpha^Q }G)\), i.e. \(\bigcap_{\mathcal{P} \in R_Q} \mathcal{P} = \mathcal{I}_Q\rtimes _\alpha G\).
	\end{enumerate}
\end{proposition}
\begin{proof}
	\(1 \Rightarrow  2.\) By Proposition \ref{prop:prim_is_fibered} the map \(\sigma \) is continuous. If \(\sigma \) is also relatively open, then for every subset \(A \subset \sigma (\operatorname{Prim}(\mathcal{A}\rtimes _{\alpha }G))\) we have \(\overline{\sigma ^{-1} (A)}  = \sigma ^{-1}(\overline{A} \cap \sigma(\operatorname{Prim}(\mathcal{A}\rtimes _{\alpha }G))) = \sigma ^{-1}(\overline{A} )\), in particular for \(A = \{ Q \} \).
	
	\(2 \Rightarrow  3.\) To prove this it suffices to observe that \(\sigma ^{-1}(Q) = R_Q \) and \(\sigma ^{-1}(\overline{\{ Q \}}) = \operatorname{Prim}(\mathcal{A}_{Q}\rtimes _{\alpha^Q }G)\).

	\(3 \Rightarrow  1.\) The proof of relative openness of \(\sigma \) generalizes the proof of \cite[Theorem~8.6]{Will2007}. Let \(U \subset \operatorname{Prim} (\mathcal{A}\rtimes _{\alpha }G)\) be an open set, and let \(\{ Q_n \} \subset \sigma (\operatorname{Prim}(\mathcal{A}\rtimes _{\alpha }G)) \setminus \sigma (U) \) be a net converging to an element \(Q_0\) in the topology of \(\sigma (\operatorname{Prim}(\mathcal{A}\rtimes _{\alpha }G))\). It is enough to show that \(Q_0 \notin \sigma (U)\). Assume the contrary, namely, that there exists a primitive ideal \(\mathcal{P}_0 \in U\) such that \(\sigma (\mathcal{P}_0) = Q_0 \). Since the set \(U\) is open in the Jacobson topology, there exists an ideal \(\mathcal{J} \in \operatorname{I}(\mathcal{A}\rtimes _\alpha G) \) such that
	\begin{equation}\label{eq:U-not-J-neighbourhood}
		U=  \{ \mathcal{P} \in \operatorname{Prim} (\mathcal{A}\rtimes _\alpha G)\colon \mathcal{P} \not\supset \mathcal{J} \}.
	\end{equation}
	We will show that for some \(n\) there exists a primitive ideal \(\mathcal{P}_n \in \operatorname{Prim} (\mathcal{A}\rtimes _\alpha G) \) such that \( \mathcal{P}_n \in R_{Q_n} \), i.e. \(\sigma (\mathcal{P}_n)=Q_n \not \in \sigma (U)\), but at the same time \(\mathcal{P}_n \not \supset \mathcal{J}\), i.e. \(\sigma (\mathcal{P}_n) \in \sigma  (U)\), and thus we obtain a contradiction.
	
	Recall that the map \(\kappa \colon
	\operatorname{Prim} \mathcal{A} \to \mathcal{Q}_{\mathcal{A}}\) is continuous and
	open. Therefore, from \(Q_0 \in \overline{\{Q_n\}}\) it follows that \(\kappa ^{-1} (Q_0) \subset
	\overline{\kappa ^{-1} (\{ Q_n \})}\).
	By the definition of the Jacobson topology on
	\(\operatorname{Prim} \mathcal{A}\), the latter holds if and only if the inclusion
	\[ 
		\bigcap_n \mathcal{I}_{Q_n}  = \bigcap_{P \in\overline{\kappa ^{-1} (\{ Q_n \})} } P \subset \bigcap_{P \in{\kappa ^{-1} (Q_0)} } P = \mathcal{I}_{Q_0}
	\]
	holds. By amenability of the group \(G\) and the properties of the map \(\operatorname{Ind} \), listed before the statement of the proposition, it follows that
	\[
		\bigcap_{n} (\mathcal{I}_{Q_n} \rtimes _\alpha G )\subset \mathcal{I}_{Q_0} \rtimes _\alpha G. 
	\]
	From this and from the relations \( \mathcal{I}_{Q_0} \rtimes _\alpha G \subset \mathcal{P}_0\in U\) and \eqref{eq:U-not-J-neighbourhood} it follows that for some \(n\) we have \(\mathcal{I}_{Q_n} \rtimes _\alpha G \not \supset \mathcal{J}\). By condition 3 the ideal \(\mathcal{I}_{Q_n} \rtimes _\alpha G \) coincides with the intersection of all primitive ideals from the set \(R_{Q_n} \), hence there exists an ideal \(\mathcal{P}_n \in \operatorname{Prim} (\mathcal{A} \rtimes _\alpha G)\) from \(R_{Q_n}\) such that \(\mathcal{P}_n \not \supset \mathcal{J}\), as required.
\end{proof}

\begin{corollary} \label{cor:localization-open-cont-surj}
	Under the assumptions of Theorem \ref{th:orbit-localisation}, the map \(\sigma \colon \operatorname{Prim} (\mathcal{A}\rtimes _{\alpha }G) \to \mathcal{Q}_{\mathcal{A}}\) is continuous and surjective. If, in addition, the group \(G\) is amenable, then \(\sigma\) is open.
\end{corollary}
\begin{proof}
	The map \(\sigma\) was defined only in the case of quasi-regularity of the dynamical system \((\mathcal{A},G,\alpha )\). By Proposition \ref{prop:prim_is_fibered}, the map \(\sigma\) is always continuous. The condition that each quasi-orbit is locally closed implies the surjectivity of \(\sigma\) (Lemma \ref{lm:sigma-surjective}). Let us check that condition 3 of Proposition \ref{prop:prim_is_fibered:opennes} holds; then, the openness of \(\sigma\) follows from its surjectivity.
	
	Let \(Q \in \sigma (\operatorname{Prim} (\mathcal{A}\rtimes _\alpha G))\).
	Since the primitive spectrum \(\operatorname{Prim} (\tilde{\mathcal{J}}_Q) = \kappa ^{-1}(Q)\) is dense in the primitive spectrum \(\operatorname{Prim} (\mathcal{A}_Q) = \overline{\kappa ^{-1}(Q)}\), it follows from Lemma \ref{lm:essentiality-dense} (see below) that the ideal \(\tilde{\mathcal{J}}_Q\) in \(\mathcal{A}_Q\) is essential. 
	Using Lemma \ref{lm:J-essential-JxG-essential}, we then obtain that the ideal \(\tilde{\mathcal{J}}_Q \rtimes _{\alpha^Q} G\) in \(\mathcal{A}_{Q}\rtimes _{\alpha^Q }G\) is essential, which by Lemma \ref{lm:essentiality-dense} is equivalent to the density of \(\operatorname{Prim} (\tilde{\mathcal{J}}_Q \rtimes _{\alpha^Q} G) = R_Q\) in the space \(\operatorname{Prim} (\mathcal{A}_{Q}\rtimes _{\alpha^Q }G)\).
\end{proof}
Recall that an ideal \(\mathcal{J}\) in an algebra \(\mathcal{A}\) is called essential if for every nontrivial ideal \(\mathcal{I} \subset \mathcal{A}\) one has \(\mathcal{I} \cap \mathcal{J} \neq \{ 0 \}\).
\begin{lemma}	\label{lm:essentiality-dense}
	A closed two-sided ideal \(\mathcal{J}\) in a \(C^*\)-algebra \(\mathcal{A}\) is essential if and only if the subset \(\operatorname{Prim} \mathcal{J}\) is dense in the space \(\operatorname{Prim} \mathcal{A}\).
\end{lemma}
\begin{proof}
	Recall that \(\operatorname{Prim} \mathcal{J} = \{P \in \operatorname{Prim} \mathcal{A}: P \not \supset \mathcal{J} \}\). The equality \(\bigcap_{P \in \operatorname{Prim} \mathcal{J}}P \cap  \mathcal{J} = \{ 0 \}\) holds. If the ideal \(\mathcal{J}\) is essential, it follows that \(\bigcap_{P \in \operatorname{Prim} \mathcal{J}}P = \{ 0 \}\), which is equivalent to the density of \(\operatorname{Prim} \mathcal{J}\) in \(\operatorname{Prim} \mathcal{A}\).  

	Conversely, suppose that \(\bigcap_{P \in \operatorname{Prim} \mathcal{J}}P = \{ 0 \}\). Assume that for some ideal \(\mathcal{I}\) one has \(\mathcal{I} \cap \mathcal{J} = \{ 0 \}\). 
	Note that if \(P \in \operatorname{Prim} \mathcal{J}\), then such an ideal \(\mathcal{I}\) is contained in \(P\). Indeed, every primitive ideal is prime (see \cite[II.6.1.11]{Blac2006}), therefore from the relations \(\{ 0 \} =\mathcal{I} \cap \mathcal{J} \subset P\) it follows (by the definition of prime ideals) that either \(\mathcal{I} \subset P\), or \(\mathcal{J} \subset P\), but the latter contradicts the choice of \(P\). Thus, \(\mathcal{I} \subset \bigcap_{P \in \operatorname{Prim} \mathcal{J}}P = \{ 0 \}\), and hence \(\mathcal{J}\) is an essential ideal.
\end{proof}
\begin{lemma}\label{lm:J-essential-JxG-essential}
	Let \((\mathcal{A},G,\alpha )\) be a \(C^*\)-dynamical system with \(G\) amenable. If \(\mathcal{J}\) is an essential invariant ideal in \(\mathcal{A}\), then the ideal \(\mathcal{J} \rtimes _\alpha G\) in \(\mathcal{A}\rtimes _\alpha G\) is also essential.
\end{lemma}
\begin{proof}
	Let \(\pi'\) be some nondegenerate faithful representation of the ideal \(\mathcal{J}\) on a Hilbert space \(\mathcal{H}\).
	There exists a unique extension of the representation \(\pi'\) to a nondegenerate representation of the algebra \(\mathcal{A}\) on the same Hilbert space; denote it by \(\pi\). Since \(\mathcal{J}\) is an essential ideal, from the equalities \((\ker \pi) \cap \mathcal{J}= \ker \pi ' = \{ 0 \} \) it follows that \(\ker \pi =\{ 0 \}\), hence the representation \(\pi\) is faithful. Let \(\tilde{\pi}\rtimes \lambda\) be the corresponding representation of the algebra \(\mathcal{A} \rtimes _\alpha G\), introduced before Proposition \ref{prop:prim_is_fibered:opennes}.
	Since the group \(G\) is amenable, we have equalities \(\ker (\tilde{\pi}\rtimes \lambda) = \operatorname{Ind} (\ker \pi) = \{ 0 \} \rtimes _\alpha G = \{ 0 \}\). In other words, the representation \(\tilde{\pi}\rtimes \lambda\) is faithful. 
	The representation \(\pi\) is nondegenerate, therefore the representation \(\tilde{\pi}\rtimes \lambda\) is also nondegenerate (\cite[Lemma~2.17]{Will2007}). For the same reason, the representation \(\widetilde{\pi '} \rtimes \lambda\) is nondegenerate as well. 
	Moreover, the restriction of \(\tilde{\pi}\rtimes \lambda\) to the ideal \(\mathcal{J} \rtimes _\alpha G\) is faithful and coincides with the representation \(\widetilde{\pi '} \rtimes \lambda\). Thus, we may regard the algebra \( \mathcal{A} \rtimes _\alpha G\) as a nondegenerate subalgebra of the algebra of bounded operators \(B(L_2(G,\mathcal{H}))\), containing \(\mathcal{J} \rtimes _\alpha G\) as an ideal, which is also a nondegenerate subalgebra of operators. We now show that the ideal \(\mathcal{J}\rtimes _\alpha G\) is essential. Let \(\mathcal{I}\) be an ideal in \( \mathcal{A} \rtimes _\alpha G\) such that \(\mathcal{I} \cap (\mathcal{J} \rtimes _\alpha G) = \{ 0 \}\). Then, from the equality \((\mathcal{J} \rtimes _\alpha G) \cdot \mathcal{H}' = \mathcal{H}'\), where \(\mathcal{H}' = L_2(G,\mathcal{H})\), we have
	\[\{ 0 \} =(\mathcal{I} \cap (\mathcal{J} \rtimes _\alpha G)) \cdot \mathcal{H}' = \mathcal{I} \cdot  (\mathcal{J} \rtimes _\alpha G)\cdot \mathcal{H}' = \mathcal{I} \cdot  ((\mathcal{J} \rtimes _\alpha G) \cdot \mathcal{H}') = \mathcal{I} \cdot \mathcal{H}', \]
	that is, \(\mathcal{I} = \{ 0 \}\).	
\end{proof}

\section{Free action of a group on the spectrum}
\label{sec:action-free}

The action of a group \(G\) on the primitive spectrum \(\operatorname{Prim} \mathcal{A}\) is called free if for all \(P \in \operatorname{Prim} \mathcal{A}\) the stabilizer \(G_P=\{ t\in G \colon \alpha _t(P)=P \} \) is trivial. In this section we show that the map \(\sigma\), studied in the previous section, is a homeomorphism if the action of the group \(G\) is free on the primitive spectrum \(\operatorname{Prim} \mathcal{A}\) and certain additional conditions are satisfied.

Recall that a dynamical system \((\mathcal{A},G,\alpha )\) is called regular if it is quasi-regular, every quasi-orbit \(Q \in \mathcal{Q}_{\mathcal{A}}\) is locally closed, and for all \(Q \in \mathcal{Q}_{\mathcal{A}}\) there exists a primitive ideal \(P \in \kappa ^{-1}(Q)\) such that the map \(G/G_P \to \kappa ^{-1}(Q)\colon t\cdot G_P \mapsto \alpha _t(P)\) is a homeomorphism (in other words, the quasi-orbit consists of a single orbit homeomorphic to \(G/G_P\)).

\begin{proposition}\label{prop:free-action-regular-orbit}
	Let \((\mathcal{A},G,\alpha )\) be a regular \(C^*\)-dynamical system, and suppose the action of the group \(G\) on the primitive spectrum \(\operatorname{Prim} \mathcal{A}\) is free.
	Then, the map \(\sigma \colon \operatorname{Prim} (\mathcal{A}\rtimes _{\alpha }G) \to \mathcal{Q}_{\mathcal{A}}\) is a homeomorphism, and the inverse map to \(\sigma\) is given by \(Q \mapsto \operatorname{Ind} \mathcal{I}_Q\).
\end{proposition}
\begin{proof}
	It follows from \cite[Proposition~20]{Gree1978} that under the assumptions of the present proposition every primitive ideal \(\mathcal{P} \in \operatorname{Prim} (\mathcal{A}\rtimes _\alpha G)\) has the form \(\operatorname{Ind} P\), where \(P \in \operatorname{Prim} \mathcal{A}\).
	Since the map \(\operatorname{Ind}\) sends intersections of ideals to intersections, and every ideal is an intersection of primitive ideals, the map \(\operatorname{Ind} \colon \operatorname{I} (\mathcal{A}) \to \operatorname{I}(\mathcal{A}\rtimes _\alpha G)\) is surjective. 
	Hence, by \cite[Proposition~11]{Gree1978}, the map \(\operatorname{Ind}\) restricts to a homeomorphism between the subspace \(\operatorname{I}^G(\mathcal{A})\) of invariant ideals in \(\mathcal{A}\) and the space \(\operatorname{I} (\mathcal{A}\rtimes _\alpha G)\).
	Moreover, for every ideal \(\mathcal{I} \in \operatorname{I}^G(\mathcal{A})\) one has the equality \(\operatorname{Res} \operatorname{Ind} \mathcal{I} = \mathcal{I}\), where \(\operatorname{Res}\) is the map defined in the proof of Proposition~\ref{prop:prim_is_fibered}. 
	Thus, the continuous map \(\sigma = q \circ \operatorname{Res}\) is injective, and its inverse is given by \(\sigma (\operatorname{Prim} (\mathcal{A}\rtimes _\alpha G)) \ni Q \mapsto \operatorname{Ind} \mathcal{I}_Q\), which is also continuous.
	By the definition of regularity of a dynamical system and Lemma~\ref{lm:sigma-surjective}, the map \(\sigma\) is surjective, hence it is a homeomorphism.
\end{proof}

\begin{remark}
	Let a \(C^*\)-dynamical system \((\mathcal{A},G,\alpha)\) be quasi-regular, let the group \(G\) be compact, let its action on the primitive spectrum \(\operatorname{Prim} \mathcal{A}\) be free, and suppose that for every quasi-orbit \(Q \in \mathcal{Q}_{\mathcal{A}}\) there exists a primitive ideal \(P \in \kappa^{-1}(Q)\) which is locally closed in \(\operatorname{Prim} \mathcal{A}\).
	Then, by \cite[Lemma~22]{Gree1978}, the dynamical system \((\mathcal{A},G,\alpha)\) is regular.
\end{remark}

A dynamical system \((\mathcal{A},G,\alpha)\) is called separable if the algebra \(\mathcal{A}\) is separable and the group \(G\) is second countable. By Lemma~\ref{lm:quasi-regularity}, a separable dynamical system is quasi-regular. 
\begin{remark}
	If a dynamical system \((\mathcal{A},G,\alpha)\) is separable, then Proposition~\ref{prop:free-action-regular-orbit} reduces to \cite[Theorem~24]{Gree1978} (condition~1), whose proof is partially reproduced above. Moreover, \cite[Corollary~7.35]{Echt2017} contains a statement analogous to Proposition~\ref{prop:free-action-regular-orbit} for smooth dynamical systems, which form a slightly narrower class of dynamical systems than regular ones. However, the proof of \cite[Corollary~7.35]{Echt2017} for smooth dynamical systems admits a generalization to regular dynamical systems. Indeed, the proof of \cite[Corollary~7.35]{Echt2017} relies on \cite[Proposition~7.26]{Echt2017}, which also holds for regular dynamical systems.
\end{remark}
The following statement is also contained in \cite[Corollary~7.35]{Echt2017}.
\begin{proposition} \label{prop:free-action-separable}
	Let \((\mathcal{A},G,\alpha)\) be a separable \(C^*\)-dynamical system, let the group \(G\) be amenable, and let its action on the primitive spectrum \(\operatorname{Prim} \mathcal{A}\) be free.
	Then,  the map \(\sigma \colon \operatorname{Prim} (\mathcal{A}\rtimes _{\alpha} G) \to \mathcal{Q}_{\mathcal{A}}\) is a homeomorphism, and the inverse map to \(\sigma\) is given by \(Q \mapsto \operatorname{Ind} \mathcal{I}_Q = \mathcal{I}_Q \rtimes _\alpha G\).
\end{proposition}
\begin{proof}
	By \cite[Corollary~3.2]{Goot.Rose1979}, every ideal \(\mathcal{P} \in \operatorname{Prim} (\mathcal{A}\rtimes _\alpha G)\) has the form \(\operatorname{Ind} P\), where \(P \in \operatorname{Prim} \mathcal{A}\). 
	Arguing as in the proof of Proposition~\ref{prop:free-action-regular-orbit}, we conclude that the continuous map \(\sigma\) is injective, and the inverse map is given by 
	\(\sigma (\operatorname{Prim} (\mathcal{A}\rtimes _\alpha G)) \ni Q \mapsto \operatorname{Ind} \mathcal{I}_Q\) and is also continuous. 
	Repeating the argument from the second paragraph of the proof of \cite[Theorem~24]{Gree1978}, one can show that for every \(Q\), the ideal \(\operatorname{Ind} \mathcal{I}_Q\) is primitive.
\end{proof}

The action of a group \(G\) on the primitive spectrum \(\operatorname{Prim} \mathcal{A}\) is called essentially free if for every closed invariant subspace \(Y \subset \operatorname{Prim} \mathcal{A}\), the set \(\{ P \in Y : G_P = \{ e \} \}\) is dense in \(Y\).

\begin{proposition} \label{prop:free-action-discrete}
	Let \((\mathcal{A},G,\alpha)\) be a \(C^*\)-dynamical system, let the group \(G\) be discrete and amenable, let its action be essentially free on the primitive spectrum \(\operatorname{Prim} \mathcal{A}\), and let the space \(\mathcal{Q}_{\mathcal{A}}\) be second countable or almost Hausdorff. 
	Then, the \(C^*\)-dynamical system \((\mathcal{A},G,\alpha)\) is quasi-regular, the map \(\sigma \colon \operatorname{Prim} (\mathcal{A}\rtimes _{\alpha} G) \to \mathcal{Q}_{\mathcal{A}}\) is a homeomorphism, and the inverse map to \(\sigma\) is given by \(Q \mapsto \operatorname{Ind} \mathcal{I}_Q = \mathcal{I}_Q \rtimes _\alpha G\).
\end{proposition}
\begin{proof}
The proof is obtained by a literal repetition of the proof of \cite[Lemma~9]{Bayb.Sara2024}, which generalizes the proof of \cite[Corollary~2.6]{Echt.Laca2013}, established for separable dynamical systems. The mentioned \cite[Lemma~9]{Bayb.Sara2024} is stated under the assumption that the group \(G\) is abelian, but this condition can be relaxed to amenability without any obstacles. Moreover, it assumes that the quasi-orbit space is almost Hausdorff, but its proof also applies when the quasi-orbit space is second countable (while the dynamical system itself may not be separable), since in this case every irreducible closed subset of the quasi-orbit space is the closure of a point (see \cite[Lemma on p.~222]{Gree1978}).
\end{proof}

We emphasize that under the assumptions of Theorem~\ref{th:orbit-localisation}, Corollary~\ref{cor:localization-open-cont-surj}, and under the conditions of one of Propositions~\ref{prop:free-action-regular-orbit}, \ref{prop:free-action-separable}, or \ref{prop:free-action-discrete}, each set \(R_Q\) is a singleton, i.e., injectivity of the map \(\sigma\) is added to the conclusions of Corollary~\ref{cor:localization-open-cont-surj}. Clearly, \(R_Q\) will be a singleton if the dynamical system \((\tilde{\mathcal{J}}_Q, G, \alpha^Q)\) satisfies the assumptions of one of Propositions~\ref{prop:free-action-regular-orbit}, \ref{prop:free-action-separable}, or \ref{prop:free-action-discrete}. The resulting homeomorphism \(\sigma \colon \operatorname{Prim} (\mathcal{A}\rtimes _{\alpha} G) \to \mathcal{Q}_{\mathcal{A}}\) makes it possible, starting from the Jacobson topology on \(\operatorname{Prim} \mathcal{A}\), to describe the Jacobson topology on \(\operatorname{Prim} (\mathcal{A}\rtimes _{\alpha} G)\).

Note also that the results of this section are obtained in a somewhat different setting than the assumptions of Corollary~\ref{cor:localization-open-cont-surj}. In particular, Proposition~\ref{prop:free-action-regular-orbit} did not require the amenability of the group \(G\). Moreover, in Proposition~\ref{prop:free-action-separable}, instead of local closedness of quasi-orbits, separability of the dynamical system was required, while in Proposition~\ref{prop:free-action-discrete} it was assumed that the quasi-orbit space \(\mathcal{Q}_{\mathcal{A}}\) is second countable; both of these conditions are unrelated to the assumption of local closedness of quasi-orbits.

\section{An example of a dynamical system with trivial group action on the spectrum}
\label{sec:action-fixed}

We consider an example of a dynamical system in which the group action on the primitive spectrum is trivial. Let \((\mathcal{D},G,\beta )\) be a dynamical system, and let \(\mathcal{C}\) be a nuclear \(C^*\)-algebra. Then, \((\mathcal{C}\otimes \mathcal{D},G,\operatorname{id}  \otimes \beta)\) is also a dynamical system (see \cite[Remark~2.74]{Will2007}).
\begin{lemma}\label{lm:C-ox-D-fixed}
	Let \((\mathcal{C}\otimes \mathcal{D},G,\operatorname{id}  \otimes \beta)\) be a \(C^*\)-dynamical system, where \(\mathcal{C}\) is a nuclear \(C^*\)-algebra. If the \(C^*\)-algebra \(\mathcal{D}\) is simple (i.e., has no nontrivial closed two-sided ideals), then the action \(\operatorname{id}  \otimes \beta\) is trivial on \(\operatorname{Prim} (\mathcal{C}\otimes \mathcal{D})\).
\end{lemma}
\begin{proof}
	Denote by \(\operatorname{Prime} \mathcal{A}\) the subspace in \(\operatorname{I}(\mathcal{A})\) consisting of prime closed two-sided ideals of the algebra \(\mathcal{A}\).
	Since \(\mathcal{C}\) is a nuclear algebra, the map
	\begin{equation}\label{eq:Pi-Prime-tensor-prod}
	\Pi \colon  \operatorname{Prime} \mathcal{C} \times \operatorname{Prime} \mathcal{D} \to \operatorname{Prime} (\mathcal{C}\otimes \mathcal{D}) \colon (P,K)\mapsto P\otimes \mathcal{D}+ \mathcal{C}\otimes K
	\end{equation}
	is a homeomorphism (see \cite[IV.3.4.26]{Blac2006} and \cite[Proposition~2.16]{Blan.Kirc2004}). By assumption, the set \(\operatorname{Prime} \mathcal{D}\) contains only the trivial ideal.
	Hence, every prime ideal in the algebra \(\mathcal{C}\otimes \mathcal{D}\) has the form \(P\otimes \mathcal{D}\) for some \(P \in \operatorname{Prime} \mathcal{C}\). It remains to note that for every prime (in particular, primitive) ideal \(P\otimes \mathcal{D} \in  \operatorname{Prime} (\mathcal{C}\otimes \mathcal{D})\), for all \(t \in G\) we have \((\operatorname{id} \otimes \beta_t) (P \otimes \mathcal{D}) = P \otimes \mathcal{D} \), as required.
\end{proof}

\begin{corollary}
	\label{cor:A-iso-C-ox-D-fixed}
	Let \((\mathcal{A},G,\alpha )\) be a \(C^*\)-dynamical system. Suppose that there exists a \(C^*\)-dynamical system \((\mathcal{D},G,\beta )\), where \(\mathcal{D}\) is a simple \(C^*\)-algebra, and
	an isomorphism \(\phi \colon \mathcal{A} \to \mathcal{C}\otimes \mathcal{D}\), where \(\mathcal{C}\) is a nuclear \(C^*\)-algebra, such that for all \(t \in G\) the equality \(\phi \circ \alpha_t \circ \phi ^{-1}   = \operatorname{id}  \otimes \beta_t \) holds. Then, the action \(\alpha \) is trivial on \(\operatorname{Prim} \mathcal{A}\).
\end{corollary}

The map \(\Pi\) from \eqref{eq:Pi-Prime-tensor-prod} for an arbitrary algebra \(\mathcal{D}\) restricts to a continuous embedding
\begin{equation}\label{eq:Pi1-Prim-tensor-prod}
	\Pi _1 \colon  \operatorname{Prim} \mathcal{C} \times \operatorname{Prim} \mathcal{D} \to \operatorname{Prim} (\mathcal{C}\otimes \mathcal{D}) \colon (P,K)\mapsto P\otimes \mathcal{D}+ \mathcal{C}\otimes K,
\end{equation}
which sends a pair consisting of the kernels of irreducible representations \([\pi _{\mathcal{C}}] \in \widehat{\mathcal{C}}\) and \([\pi _{\mathcal{D}}] \in \widehat{\mathcal{D}}\) to the kernel of the irreducible representation \(\pi_{\mathcal{C}} \otimes \pi _{\mathcal{D}}\). The map \(\Pi _1\) is, in general, not surjective. It becomes surjective (and thus a homeomorphism), for example, when the algebras \(\mathcal{C}\) and \(\mathcal{D}\) are separable, or when one of them is of type I (see \cite[IV.3.4.22--IV.3.4.26]{Blac2006}). In what follows, we will assume that \(\Pi _1\) is a homeomorphism.

By \cite[Lemma~2.75]{Will2007}, there is a canonical isomorphism
\begin{equation}\label{eq:Pi2-Prim-tensor-prod}
	\psi\colon (\mathcal{C}\otimes \mathcal{D})\rtimes _{\operatorname{id}  \otimes \beta} G \to \mathcal{C}\otimes (\mathcal{D}\rtimes _{\beta} G).
\end{equation}
Consider the analogue of the map \eqref{eq:Pi1-Prim-tensor-prod} for the algebra \((\mathcal{C}\otimes \mathcal{D})\rtimes _{\operatorname{id}  \otimes \beta} G\), corresponding to the dynamical system from Lemma \ref{lm:C-ox-D-fixed}:
\begin{equation}
\begin{aligned}
	\Pi _2\colon \operatorname{Prim} \mathcal{C} \times \operatorname{Prim} (\mathcal{D}\rtimes _{\beta} G) &\to \operatorname{Prim} (\mathcal{C}\otimes (\mathcal{D}\rtimes _{\beta} G))\colon \\
	(P,K) &\mapsto P\otimes (\mathcal{D}\rtimes _{\beta} G)+ \mathcal{C}\otimes K.
\end{aligned}
\end{equation}
If \(\Pi _2\) is a homeomorphism, then we obtain a description of the spectrum \(\operatorname{Prim}((\mathcal{C}\otimes \mathcal{D})\rtimes _{\operatorname{id}  \otimes \beta} G)\) in terms of \(\operatorname{Prim} \mathcal{C}\) and \(\operatorname{Prim}  (\mathcal{D}\rtimes _{\beta} G)\).

Let the dynamical system \((\mathcal{A},G,\alpha )\) satisfy the assumptions of Corollary \ref{cor:A-iso-C-ox-D-fixed}. Then, the algebras \(\mathcal{A}\rtimes _\alpha G\) and \((\mathcal{C}\otimes \mathcal{D})\rtimes _{\operatorname{id}  \otimes \beta} G\) are isomorphic (see \cite[Lemma~2.65]{Will2007}); denote the corresponding isomorphism by \(\phi \rtimes \operatorname{id}\). The following obvious lemma gives a description of the primitive spectrum \(\operatorname{Prim} (\mathcal{A} \rtimes _{\alpha } G)\).
\begin{lemma}\label{lm:C-ox-D-spectrum}
	Let \((\mathcal{A},G,\alpha )\) satisfy the assumptions of Corollary \ref{cor:A-iso-C-ox-D-fixed}, and let the embeddings
	\begin{equation}\label{eq:Pi_1-inclusion}
		\Pi _1 \colon  \operatorname{Prim} \mathcal{C}  \to \operatorname{Prim} (\mathcal{C}\otimes \mathcal{D}) \colon P \mapsto P\otimes \mathcal{D},
	\end{equation} 
	\begin{equation}\label{eq:Pi_2-inclusion}
	\begin{aligned}
		\Pi _2\colon \operatorname{Prim} \mathcal{C} \times \operatorname{Prim} (\mathcal{D}\rtimes _{\beta} G) &\to \operatorname{Prim} (\mathcal{C}\otimes (\mathcal{D}\rtimes _{\beta} G))\colon \\
		(P,K) &\mapsto P\otimes (\mathcal{D}\rtimes _{\beta} G)+ \mathcal{C}\otimes K
	\end{aligned}
	\end{equation}
	be homeomorphisms. Then, the map
	\begin{align}
		\begin{split}
					\Psi \colon &\operatorname{Prim} \mathcal{A} \times \operatorname{Prim} (\mathcal{D}\rtimes _{\beta} G) \to \operatorname{Prim} (\mathcal{A}\rtimes _\alpha G) \colon \\
					& (P,K) \mapsto (\phi \rtimes \operatorname{id})^{-1} \circ \psi ^{-1} \circ \Pi _2\left((\Pi_1 ^{-1} \circ \phi (P),K)\right)
		\end{split} 
	\end{align}
	is a homeomorphism as well.
\end{lemma}
Under the assumptions of Corollary \ref{cor:A-iso-C-ox-D-fixed} every quasi-orbit in \(\operatorname{Prim} \mathcal{A}\) consists of a single point, \(\kappa^{-1}(Q)=\{P\}=\{\mathcal{I}_Q\}\), hence the spaces \(\mathcal{Q}_{\mathcal{A}}\) and \(\operatorname{Prim}\mathcal{A}\) are homeomorphic. If the hypotheses of Lemma \ref{lm:C-ox-D-spectrum} are satisfied, then using the explicit form of the map \(\Psi\) one checks that for arbitrary irreducible representations \(\pi_{\mathcal{C}}\) and \(\pi_{\mathcal{D}}' \rtimes V\) of \(\mathcal{C}\) and \(\mathcal{D}\rtimes_\beta G\), respectively, the equality
\[
	\Psi \bigl(\phi^{-1}\bigl(\Pi_1 (\ker \pi_{\mathcal{C}}) \bigr) , \ker (\pi_{\mathcal{D}}' \rtimes V)\bigr) = \ker \Bigl( \bigl((\pi_{\mathcal{C}} \otimes \pi_{\mathcal{D}}') \circ \phi \bigr)\rtimes \bigl(1_G \otimes V \bigr)\Bigr)
\]
holds, where \(1_G\) is the trivial representation of \(G\) on the Hilbert space corresponding to \(\pi_{\mathcal{C}}\).
Consequently, for \(P = \phi^{-1}\bigl(\Pi_1 (\ker \pi_{\mathcal{C}}) \bigr)\) and \(K = \ker (\pi_{\mathcal{D}}' \rtimes V)\) we have
\begin{align}
	\operatorname{Res} (\Psi(P,K)) = & \ker\bigl((\pi_{\mathcal{C}} \otimes \pi_{\mathcal{D}}') \circ \phi \bigr) = \phi^{-1}(\ker (\pi_{\mathcal{C}} \otimes \pi_{\mathcal{D}}')) = \\
	& \phi^{-1}\bigl((\ker \pi_{\mathcal{C}} )\otimes \mathcal{D}\bigr) = \phi^{-1}\bigl(\Pi_1 (\ker \pi_{\mathcal{C}}) \bigr) = P.
\end{align}
From this and the surjectivity of \(\Pi_1\) it follows that for an arbitrary ideal \(\mathcal{P} \in \operatorname{Prim} (\mathcal{A}\rtimes _\alpha G)\) one has \(\operatorname{Res} \mathcal{P} = P = \mathcal{I}_{\kappa(P)}\), where \(\Psi^{-1}(\mathcal{P})=(P,K)\); in other words, the dynamical system \((\mathcal{A},G,\alpha)\) is quasi-regular (see \eqref{Eq:lives-on-quasi-orbit}).
The map \(\sigma\colon \operatorname{Prim} (\mathcal{A}\rtimes _\alpha G) \to \mathcal{Q}_{\mathcal{A}}\) from Sec. \ref{sec:fibers} acts by \(\sigma(\mathcal{P})=\kappa(P)\), where \(P=\operatorname{Res}\mathcal{P}\), and is surjective. The fiber \(R_Q=\sigma^{-1}(Q)\) can be identified with the space \(\operatorname{Prim}(\mathcal{D}\rtimes_{\beta} G)\). Note that although \(\mathcal{D}\) is assumed to be simple, the crossed product \(\mathcal{D}\rtimes_{\beta} G\) need not be simple in general.

\section{An action with fixed points on the spectrum}
\label{sec:action-free-and-fixed}
We now describe the primitive spectrum of the crossed product corresponding to a dynamical system of a special type in which, roughly speaking, the group action on the primitive spectrum is free except for a certain fixed-point subset.

Let \((\mathcal{A},G,\alpha)\) be a dynamical system, and let \(\mathcal{J}_0\) be an invariant ideal in \(\mathcal{A}\).
We denote by \(\sigma_{\mathcal{J}_0}\colon \operatorname{Prim} (\mathcal{J}_0 \rtimes_\alpha G) \to \mathcal{Q}_{\mathcal{J}_0}\) and \(\operatorname{\mathcal{J}_0-Ind} \colon \operatorname{I}(\mathcal{J}_0) \to \operatorname{I}(\mathcal{J}_0 \rtimes G)\) the maps \(\sigma\) and \(\operatorname{Ind}\), respectively, with \(\mathcal{A}\) replaced by \(\mathcal{J}_0\) (see Section~\ref{sec:fibers}).
If the dynamical system \((\mathcal{J}_0,G,\alpha)\) satisfies the assumptions of one of Propositions~\ref{prop:free-action-regular-orbit}, \ref{prop:free-action-separable}, or~\ref{prop:free-action-discrete}, then the map \(\sigma_{\mathcal{J}_0}\) is a homeomorphism, and its inverse has the form
\begin{equation}\label{eq:sig^-1_J0}
	\sigma^{-1}_{\mathcal{J}_0}\colon \mathcal{Q}_{\mathcal{J}_0} \to \operatorname{Prim} (\mathcal{J}_0 \rtimes_\alpha G)\colon Q \mapsto \operatorname{\mathcal{J}_0-Ind} \mathcal{I}_Q.
\end{equation}
In what follows, we will assume that the map \eqref{eq:sig^-1_J0} is a homeomorphism. Let us rewrite the map \eqref{eq:sig^-1_J0} in terms of quasi-orbits in \(\mathcal{Q}_{\mathcal{A}}\) and primitive ideals in \(\mathcal{A}\rtimes_\alpha G\).

For any ideal \(\mathcal{J}\) in an arbitrary algebra \(\mathcal{A}\), let \(\operatorname{Prim}^{\mathcal{A}} \mathcal{J} = \{ P \in \operatorname{Prim}\mathcal{A} \colon\allowbreak P \not\supset \mathcal{J} \}\) and \(\operatorname{Prim}^{\mathcal{A}} \mathcal{A}/\mathcal{J} = \{ P \in \operatorname{Prim}\mathcal{A} \colon P \supset \mathcal{J} \}\). The restriction of the equivalence relation \(\sim_q\) from \(\operatorname{Prim} \mathcal{A}\) to \(\operatorname{Prim}^{\mathcal{A}} \mathcal{J}_0\) coincides with the equivalence relation \(\sim_q\) in \(\operatorname{Prim}\mathcal{J}_0\) of the form \eqref{eq:quasi-orbit-sim}, so the quasi-orbit space \(\mathcal{Q}_{\mathcal{J}_0}\) is homeomorphic to the space \(\kappa(\operatorname{Prim}^{\mathcal{A}} \mathcal{J}_0) \subset\mathcal{Q}_{\mathcal{A}}\), where \(\kappa \colon \operatorname{Prim} \mathcal{A} \to \mathcal{Q}_{\mathcal{A}} \) is the quotient map (see the beginning of Sec. \ref{subsec:crossed-prod-spec}). However, the ideal \(\mathcal{I}_Q = \bigcap_{P \in \kappa_{\mathcal{J}_0}^{-1}(Q)} P\) appearing in \eqref{eq:sig^-1_J0}, where \(\kappa_{\mathcal{J}_0} \colon \operatorname{Prim} \mathcal{J}_0 \to \mathcal{Q}_{\mathcal{J}_0} \) is the quotient map, generally differs from the ideal \(\mathcal{I}_Q = \bigcap_{P \in \kappa^{-1}(Q)} P \) associated with the dynamical system \((\mathcal{A},G,\alpha)\).
Specifically, every primitive ideal in \(\mathcal{J}_0\) has the form \(P \cap \mathcal{J}_0\), where \(P \in \operatorname{Prim} \mathcal{A}\) is such that \(P \not \supset \mathcal{J}_0\), hence the following relations hold:
\begin{equation} \label{eq:I-Q-J0=I-Q-A-cap-J0}
	\mathcal{I}_{\kappa_{\mathcal{J}_0}(P \cap \mathcal{J}_0)} = \bigcap_{t \in G}\alpha_t (P \cap \mathcal{J}_0) =  \bigcap_{t \in G}\alpha_t (P)  \cap \mathcal{J}_0  = \mathcal{I}_{\kappa(P)} \cap \mathcal{J}_0.
\end{equation}
From this we obtain
\begin{equation}\label{eq:sig^-1_J0-circ-kappa_J0}
		\sigma ^{-1}_{\mathcal{J}_0} ( \kappa _{\mathcal{J}_0} (P \cap \mathcal{J}_0)) = \operatorname{\mathcal{J}_0 - Ind} (\mathcal{I}_{\kappa_{\mathcal{J}_0}(P \cap \mathcal{J}_0)}) = \operatorname{\mathcal{J}_0 - Ind} ( \mathcal{I}_{\kappa(P)} \cap \mathcal{J}_0).
\end{equation}
If the group \(G\) is amenable, then by the properties of the mapping \(\operatorname{Ind}\) we have
\begin{align}\label{eq:J0-ind-intersection}
	\begin{split}
		\operatorname{\mathcal{J}_0 - Ind} ( \mathcal{I}_{\kappa(P)} \cap \mathcal{J}_0) =& (\mathcal{I}_{\kappa (P)}\cap  \mathcal{J}_0 ) \rtimes _\alpha G = 
		\operatorname{Ind} (\mathcal{I}_{\kappa (P)}\cap  \mathcal{J}_0 ) = \\ 
			& (\operatorname{Ind} \mathcal{I}_{\kappa (P)}) \cap (\operatorname{Ind} \mathcal{J}_0 ) =
	(\mathcal{I}_{\kappa (P)} \rtimes _\alpha G)\cap (\mathcal{J}_0 \rtimes _\alpha G).
	\end{split} 
\end{align}
Consequently, for every \(P \in\operatorname{Prim}^\mathcal{A} \mathcal{J}_0\) we have
\begin{equation}\label{eq:sigma^-1:formula}
	\sigma ^{-1}_{\mathcal{J}_0} (\kappa _{\mathcal{J}_0} (P \cap \mathcal{J}_0)) = (\mathcal{I}_{\kappa (P)} \rtimes _\alpha G)\cap (\mathcal{J}_0 \rtimes _\alpha G),
\end{equation}
which allows us to describe primitive ideals of the algebra \(\mathcal{A} \rtimes _\alpha G\).
\begin{lemma} \label{lm:Phi-on-prim-J0}
	Let \((\mathcal{A},G,\alpha )\) be a \(C^*\)-dynamical system with an amenable group \(G\), and let \(\mathcal{J}_0\) be a closed two-sided invariant ideal in \(\mathcal{A}\) such that the mapping 
	\(\sigma ^{-1}_{\mathcal{J}_0} \colon \mathcal{Q}_{\mathcal{J}_0} \to \operatorname{Prim} (\mathcal{J}_0 \rtimes _\alpha G)\) is a homeomorphism. 
	Then, the mapping 
	\begin{equation} \label{eq:Phi-on-prim-J0}
		\operatorname{Prim}^\mathcal{A} \mathcal{J}_0 \to \operatorname{Prim}(\mathcal{A} \rtimes_\alpha G)\colon P \mapsto \mathcal{I}_{\kappa(P)} \rtimes _\alpha G
	\end{equation}
	is continuous and open, and its image coincides with \(\operatorname{Prim}^{\mathcal{A} \rtimes_\alpha G} (\mathcal{J}_0 \rtimes _\alpha G)\). In particular, the mapping 
	\[
		\kappa(\operatorname{Prim}\mathcal{J}_0) \to \operatorname{Prim}^{\mathcal{A} \rtimes_\alpha G} (\mathcal{J}_0 \rtimes _\alpha G) \colon Q \mapsto \mathcal{I}_{Q} \rtimes _\alpha G
	\]
	is a homeomorphism.
\end{lemma}
\begin{proof}
	From the assumptions on the dynamical system \((\mathcal{J}_0,G,\alpha )\), the equality \eqref{eq:sigma^-1:formula}, and the openness of the mapping \(\kappa_{\mathcal{J}_0}\), we obtain the continuity, surjectivity, and openness of the mapping 
	\begin{equation}\label{Phi-on-prim-j0:2}
		\sigma ^{-1}_{\mathcal{J}_0} \circ \kappa _{\mathcal{J}_0} \colon \operatorname{Prim} \mathcal{J}_0 \to \operatorname{Prim} (\mathcal{J}_0 \rtimes _\alpha G)\colon P \cap \mathcal{J}_0 \mapsto (\mathcal{I}_{\kappa (P)} \rtimes _\alpha G)\cap (\mathcal{J}_0 \rtimes _\alpha G).
	\end{equation}
	Since the correspondences 
	\begin{align}
		\operatorname{Prim}^\mathcal{A} \mathcal{J}_0 \to & \operatorname{Prim} \mathcal{J}_0\colon P \mapsto P \cap \mathcal{J}_0, \\
		\operatorname{Prim}^{\mathcal{A} \rtimes_\alpha G} (\mathcal{J}_0 \rtimes _\alpha G) \to & \operatorname{Prim} (\mathcal{J}_0 \rtimes _\alpha G) \colon  \mathcal{P} \mapsto \mathcal{P} \cap (\mathcal{J}_0 \rtimes _\alpha G)
	\end{align} 
	are homeomorphisms, to prove the lemma it suffices to show that the ideal \(\mathcal{I}_{\kappa (P)} \rtimes _\alpha G\) is primitive in the algebra \(\mathcal{A} \rtimes_{\alpha} G\) for all \(P \in\operatorname{Prim}^\mathcal{A} \mathcal{J}_0 \).
	
	Instead of proving the primitivity of the ideal \(\mathcal{I}_{\kappa (P)} \rtimes _\alpha G\) in the algebra \(\mathcal{A} \rtimes _\alpha G\), we will check the primitivity of the ideal \(\{ 0 \}\) in the quotient algebra \(\mathcal{A}_{\kappa (P)} \rtimes _{\alpha^{\kappa (P)}} G\).
	Let \(\mathcal{J}_{0,\kappa (P)}\) be the image of the ideal \(\mathcal{J}_0\) in the quotient algebra \(\mathcal{A}_{\kappa (P)}\).
	From the equality \eqref{eq:sigma^-1:formula} we obtain that the ideal \((\mathcal{I}_{\kappa (P)} \rtimes _\alpha G)\cap (\mathcal{J}_0 \rtimes _\alpha G)\) is primitive in \(\mathcal{J}_0 \rtimes _\alpha G\), and hence \(\{ 0 \}\) is primitive in \(\mathcal{J}_{0,\kappa (P)} \rtimes _{\alpha^{\kappa (P)}} G\). We will show that \(\mathcal{J}_{0,\kappa (P)} \rtimes _{\alpha^{\kappa (P)}} G\) is an essential ideal in \(\mathcal{A}_{\kappa (P)} \rtimes _{\alpha^{\kappa (P)}} G\), which will imply that the ideal \(\{0 \}\) is primitive in \(\mathcal{A}_{\kappa (P)} \rtimes _{\alpha^{\kappa (P)}} G\).
	To this end, note that the relations
	\begin{equation}
	\kappa^{-1} (\kappa(P)) \subset  (\operatorname{Prim}^{\mathcal{A}} \mathcal{J}_0) \cap (\operatorname{Prim}^{\mathcal{A}} \mathcal{A}_{\kappa (P)})
	\end{equation} 
	and \(\operatorname{Prim}^{\mathcal{A}} \mathcal{A}_{\kappa (P)} = \overline{\kappa^{-1}(\kappa (P))}\)
	hold.
	Since \( (\operatorname{Prim}^{\mathcal{A}} \mathcal{J}_0) \cap (\operatorname{Prim}^{\mathcal{A}} \mathcal{A}_{\kappa (P)})\) and \(\operatorname{Prim}^{\mathcal{A}_{\kappa(P)}} \mathcal{J}_{0,\kappa (P)}\) are homeomorphic, the ideal \(\mathcal{J}_{0,\kappa (P)}\) is nontrivial and its spectrum is dense in \(\operatorname{Prim} \mathcal{A}_{\kappa (P)}\). It then follows from Lemmas \ref{lm:essentiality-dense} and \ref{lm:J-essential-JxG-essential} that \(\mathcal{J}_{0,\kappa (P)} \rtimes _{\alpha^{\kappa (P)}} G\) is an essential ideal in \(\mathcal{A}_{\kappa (P)} \rtimes _{\alpha^{\kappa (P)}} G\), as required.
\end{proof}
Thus, under the assumptions of Lemma \ref{lm:Phi-on-prim-J0} we obtain a map describing a part of the spectrum of the algebra \(\mathcal{A} \rtimes _\alpha G\). If, in addition, the dynamical system \((\mathcal{A}/\mathcal{J}_0, G, \alpha ^{\mathcal{J}_0})\) satisfies the hypotheses of Corollary \ref{cor:A-iso-C-ox-D-fixed} and Lemma \ref{lm:C-ox-D-spectrum}, then using the homeomorphism from Lemma \ref{lm:C-ox-D-spectrum} and the homeomorphisms
\begin{gather}
	\label{eq:quotient-J0}
	\operatorname{Prim}^{\mathcal{A}} \mathcal{A}/\mathcal{J}_0 \to \operatorname{Prim}\mathcal{A}/\mathcal{J}_0 \colon P \mapsto \rho_{\mathcal{J}_0} (P), \\
	\label{eq:quotient-J0xG}
	\begin{gathered}
		\operatorname{Prim}((\mathcal{A}/\mathcal{J}_0) \rtimes _{\alpha^{\mathcal{J}_0}} G) \to \operatorname{Prim}^{\mathcal{A}\rtimes_\alpha G} ((\mathcal{A}/\mathcal{J}_0) \rtimes _{\alpha^{\mathcal{J}_0}} G) \colon \\
		\mathcal{P} \mapsto \rho _{\mathcal{J}_0 \rtimes _\alpha G}^{-1} (\mathcal{P}),
	\end{gathered}
\end{gather} 
where \(\rho _{\mathcal{J}_0} \colon \mathcal{A} \to \mathcal{A}/\mathcal{J}_0\) and \(\rho _{\mathcal{J}_0 \rtimes _\alpha G} \colon \mathcal{A} \rtimes_\alpha G \to (\mathcal{A}/\mathcal{J}_0) \rtimes _{\alpha^{\mathcal{J}_0}} G\) are the quotient maps,
we obtain the remaining part \(\operatorname{Prim}^{\mathcal{A} \rtimes _\alpha G} ((\mathcal{A}/\mathcal{J}_0 )\rtimes _{\alpha^{\mathcal{J}_0}} G)\). The following statement allows one to glue these two maps together.
\begin{theorem}\label{th:free-and-fixed-action}
	Let \((\mathcal{A},G,\alpha )\) be a \(C^*\)-dynamical system with an amenable group \(G\), and let \(\mathcal{J}_0\) be a closed two-sided invariant ideal in \(\mathcal{A}\) such that the map 
	\(\sigma ^{-1}_{\mathcal{J}_0} \colon \mathcal{Q}_{\mathcal{J}_0} \to \operatorname{Prim} (\mathcal{J}_0 \rtimes _\alpha G)\) is a homeomorphism, and suppose the \(C^*\)-dynamical system \((\mathcal{A}/\mathcal{J}_0, G, \alpha ^{\mathcal{J}_0})\) satisfies the hypotheses of Corollary \ref{cor:A-iso-C-ox-D-fixed} and Lemma \ref{lm:C-ox-D-spectrum} for some \(\mathcal{C}\), \(\mathcal{D}\) and \(\beta\). 
	In addition, let
	\(\Phi \colon \operatorname{Prim} \mathcal{A} \times \operatorname{Prim} (\mathcal{D}\rtimes _{\beta} G) \to \operatorname{Prim} (\mathcal{A}\rtimes _\alpha G)\)
	be the map defined by
	\begin{align*}
		& \Phi (P,K)  = \operatorname{Ind} \mathcal{I}_{\kappa  (P)}   = \mathcal{I}_{\kappa (P)} \rtimes _\alpha G, & P \in & \operatorname{Prim}^{\mathcal{A}} \mathcal{J}_0, \\
		& \Phi (P,K)  = \rho _{\mathcal{J}_0 \rtimes _\alpha G}^{-1} (\Psi (\rho _{\mathcal{J}_0}( P),K)),  & P \in &\operatorname{Prim}^{\mathcal{A}} \mathcal{A}/\mathcal{J}_0 ,  K \in \operatorname{Prim} (\mathcal{D}\rtimes _{\beta} G),
	\end{align*}
	where \(\Psi \) is the map from Lemma \ref{lm:C-ox-D-spectrum}. Then, \(\Phi \) is a continuous open surjection.
\end{theorem}

\begin{proof}
	The map \(\Phi\) is surjective. Indeed, by Lemma \ref{lm:Phi-on-prim-J0} we have \(\Phi (\operatorname{Prim}^{\mathcal{A}} \mathcal{J}_0 \times  \operatorname{Prim} (\mathcal{D}\rtimes _{\beta} G)) = \operatorname{Prim}^{\mathcal{A}\rtimes_\alpha G} (\mathcal{J}_0 \rtimes _\alpha G)\), and since the maps \eqref{eq:quotient-J0} and \eqref{eq:quotient-J0xG} are homeomorphisms, by Lemma \ref{lm:C-ox-D-spectrum} the restriction of \(\Phi\) to the closed subspace \(\operatorname{Prim}^{\mathcal{A}} (\mathcal{A}/\mathcal{J}_0) \times \operatorname{Prim} (\mathcal{D}\rtimes _{\beta} G)\) is a homeomorphism onto the closed subspace \(\operatorname{Prim}^{\mathcal{A}\rtimes_\alpha G} ((\mathcal{A}/\mathcal{J}_0) \rtimes _{\alpha^{\mathcal{J}_0}} G)\). 

	Note that by Lemma \ref{lm:Phi-on-prim-J0} the restriction of \(\Phi\) to the open subspace \(\operatorname{Prim}^{\mathcal{A}} \mathcal{J}_0 \times \operatorname{Prim} (\mathcal{D}\rtimes _{\beta} G)\) yields a continuous open surjection onto the open subspace \(\operatorname{Prim}^{\mathcal{A} \rtimes _\alpha G} (\mathcal{J}_0 \rtimes _\alpha G)\). 
	Thus, to prove openness and continuity of \(\Phi\) it suffices to study images of neighbourhoods under \(\Phi\) and \(\Phi^{-1}\) only from the closed subspaces \(\operatorname{Prim}^{\mathcal{A}} (\mathcal{A}/\mathcal{J}_0) \times \operatorname{Prim} (\mathcal{D}\rtimes _{\beta} G)\) and \(\operatorname{Prim}^{\mathcal{A} \rtimes _\alpha G} ((\mathcal{A}/\mathcal{J}_0) \rtimes _{\alpha^{\mathcal{J}_0}} G)\), respectively.

	Before proving continuity of \(\Phi\), we establish auxiliary facts: the dynamical system \((\mathcal{A},G,\alpha )\) is quasi-regular, and the map \( \sigma \colon  \operatorname{Prim} (\mathcal{A}\rtimes _{\alpha }G) \to \mathcal{Q}_{\mathcal{A}}\) from Sec.~\ref{sec:fibers} is open. Indeed, for \(P \in \operatorname{Prim}^{\mathcal{A}} \mathcal{J}_0\) we have \(\operatorname{Res} (\Phi (P,K)) = \mathcal{I}_{\kappa (P)}\), and for \(P \in \operatorname{Prim}^{\mathcal{A}} (\mathcal{A}/\mathcal{J}_0)\) and \(K \in \operatorname{Prim}(\mathcal{D} \rtimes _\beta G)\) from the equalities
	\[
		\operatorname{Res} (\Phi (P,K)) = \operatorname{Res}(\rho _{\mathcal{J}_0 \rtimes _\alpha G}^{-1} (\Psi (\rho _{\mathcal{J}_0}( P),K))) = \rho_{\mathcal{J}_0}^{-1} \operatorname{Res} (\Psi (\rho _{\mathcal{J}_0}( P),K)) 
	\]
	and the remarks after Lemma \ref{lm:C-ox-D-spectrum} we obtain \(\operatorname{Res} (\Phi (P,K)) = P\), which proves quasi-regularity. To prove openness of \(\sigma\) it suffices to verify condition 3 of Proposition \ref{prop:prim_is_fibered:opennes}.
	For \(Q \in \kappa (\operatorname{Prim}^{\mathcal{A}} \mathcal{J}_0)\) we have \(R_Q = \{\mathcal{I}_Q \rtimes _\alpha G\}\), hence condition 3 holds for such \(Q\).
	Now consider \(Q \in \kappa (\operatorname{Prim}^{\mathcal{A}} (\mathcal{A}/ \mathcal{J}_0))\). Since \(\kappa\) induces a bijection between \(\operatorname{Prim}^{\mathcal{A}} (\mathcal{A}/ \mathcal{J}_0)\) and \(\kappa (\operatorname{Prim}^{\mathcal{A}} (\mathcal{A}/ \mathcal{J}_0))\), we have \(\{\mathcal{I}_Q\} = \kappa ^{-1}(Q) \). 
	The restriction of \(\Phi\) to the closed subspace \(\operatorname{Prim}^{\mathcal{A}} (\mathcal{A}/\mathcal{J}_0) \times \operatorname{Prim} (\mathcal{D}\rtimes _{\beta} G)\) is a homeomorphism onto \(\operatorname{Prim}^{\mathcal{A} \rtimes _\alpha G}((\mathcal{A}/\mathcal{J}_0) \rtimes _{\alpha ^{\mathcal{J}_0} } G )\), therefore density of \(R_Q\) in \(\operatorname{Prim}^{\mathcal{A} \rtimes _\alpha G} (\mathcal{A}_{Q}\rtimes _{\alpha^Q }G)\) is equivalent to density of the set  
	\[\Phi ^{-1}(R_{Q}) = \{ \mathcal{I}_Q \} \times  \operatorname{Prim} (\mathcal{D}\rtimes _{\beta} G) \]
	in 
	\begin{multline*}
		\Phi ^{-1}(\operatorname{Prim}^{\mathcal{A} \rtimes _\alpha G} (\mathcal{A}_{Q}\rtimes _{\alpha^Q }G)) \\
		= \{ P \in \operatorname{Prim}^{\mathcal{A}} (\mathcal{A}/\mathcal{J}_0) \colon P \supset \mathcal{I}_Q\} \times  \operatorname{Prim} (\mathcal{D}\rtimes _{\beta} G),
	\end{multline*}
	which is obvious, since \(\{ P \in \operatorname{Prim}^{\mathcal{A}} (\mathcal{A}/\mathcal{J}_0) \colon  P \supset \mathcal{I}_Q \} = \overline{\{ \mathcal{I}_Q\}}\).

	Let us verify continuity of the map \(\Phi\). 
	Let \(U\) be an open neighbourhood of the point \(\Phi (P,K) \in  \operatorname{Prim}^{\mathcal{A} \rtimes _\alpha G} ((\mathcal{A}/\mathcal{J}_0) \rtimes_{\alpha^{\mathcal{J}_0}} G)\). We will exhibit an open neighbourhood \(V \subset  \operatorname{Prim} \mathcal{A} \times \operatorname{Prim}
	(\mathcal{D}\rtimes _{\beta} G)\) of the point \((P,K)\) such that \(\Phi (V) \subset U\). 
	From openness of the map \(\sigma \colon  \operatorname{Prim} (\mathcal{A}\rtimes _{\alpha }G) \to \mathcal{Q}_{\mathcal{A}}\) it follows that the set \(V_1 = \kappa ^{-1}(\sigma (U)) \subset \operatorname{Prim} \mathcal{A}\) is open, and 
	\( (P,K) \in  V_1 \times \operatorname{Prim} (\mathcal{D}\rtimes _\beta G) \). Since the restriction of \(\sigma\) to \(\operatorname{Prim}^{\mathcal{A} \rtimes _\alpha G} (\mathcal{J}_0 \rtimes _\alpha G)\) and the map \(Q \mapsto \operatorname{Ind} \mathcal{I}_Q\) are inverse to each other, we have 
	\begin{equation} \label{eq:Phi-V1}
		\Phi \bigl((V_1 \cap \operatorname{Prim}^{\mathcal{A}} \mathcal{J}_0) \times \operatorname{Prim} (\mathcal{D}\rtimes _\beta G)\bigr) = U \cap \operatorname{Prim}^{\mathcal{A} \rtimes _\alpha G} (\mathcal{J}_0 \rtimes _\alpha G).
	\end{equation}
	On the other hand, the restriction of \(\Phi\) to the subspace \(\operatorname{Prim} ^{\mathcal{A}} (\mathcal{A}/\mathcal{J}_0) \times \operatorname{Prim} (\mathcal{D}\rtimes _{\beta} G)\) is a homeomorphism onto \(\operatorname{Prim}^{\mathcal{A} \rtimes _\alpha G} ( (\mathcal{A}/\mathcal{J}_0) \rtimes _{\alpha^{\mathcal{J}_0}} G)\). Therefore, the preimage under \(\Phi\) of the open set \(U \cap \operatorname{Prim}^{\mathcal{A} \rtimes _\alpha G} ( (\mathcal{A}/\mathcal{J}_0) \rtimes _{\alpha^{\mathcal{J}_0}} G)\) is relatively open in \(\operatorname{Prim}^{\mathcal{A}} (\mathcal{A}/\mathcal{J}_0) \times \operatorname{Prim} ({\mathcal{D}\rtimes _{\beta} G})\). In other words, there exists an open neighbourhood
	\(V_2  \subset \operatorname{Prim} \mathcal{A} \times \operatorname{Prim} (\mathcal{D}\rtimes _{\beta} G)\) of \((P,K)\) such that 
	\begin{equation}\label{eq:Phi-V2}
		\Phi (V_2 \cap (\operatorname{Prim}^{\mathcal{A}} (\mathcal{A}/\mathcal{J}_0) \times \operatorname{Prim} (\mathcal{D}\rtimes _{\beta} G))) =U \cap \operatorname{Prim}^{\mathcal{A}\rtimes_\alpha G} ( (\mathcal{A}/\mathcal{J}_0) \rtimes _{\alpha^{\mathcal{J}_0}} G).
	\end{equation}
	Put \(V = (V_1 \times \operatorname{Prim} (\mathcal{D}\rtimes _{\beta} G)) \cap V_2\). From \eqref{eq:Phi-V1}--\eqref{eq:Phi-V2} it follows that \(\Phi (V) \subset U\), and \(V\) is a neighbourhood of \((P,K)\), as required.

	Now we check openness of \(\Phi\). 
	First recall that for an arbitrary algebra \(\mathcal{A}\) the sets \(\{ U_a \colon a \in \mathcal{A}\} \) given by
	\begin{equation} \label{Eq:jacobson-top_base-system}
		U_a=\{P \in \operatorname{Prim}\mathcal{A} : \rho _P (a) \neq 0 \} = \{P \in \operatorname{Prim}\mathcal{A} : a \notin P \}
	\end{equation}
	form a base of the topology on \(\operatorname{Prim}(\mathcal{A})\), and it is enough to consider the family \(\{ U_a \colon a \in \mathcal{A}, a \geq 0\}\), since \(\left\| a^*a \right\| = \left\|  a \right\|^2\). 
	Recall also that \(\mathcal{A}\) is an essential ideal in its multiplier algebra \(M(\mathcal{A})\), hence to every primitive ideal \(P \in \operatorname{Prim}(\mathcal{A})\) there corresponds a unique primitive ideal \(\tilde{P} \in \operatorname{Prim}(M(\mathcal{A}))\) such that \(\tilde{P}\cap\mathcal{A}=P\). Thus, \(\operatorname{Prim}(\mathcal{A})\) embeds naturally into the primitive spectrum \(\operatorname{Prim}(M(\mathcal{A}))\).
	Consequently, the sets \(\{ V_a \colon a \in M(\mathcal{A})\} \) given by
	\begin{equation} \label{Eq:jacobson-top_base-system-multiplier}
	{V}_a=\{P \in \operatorname{Prim} \mathcal{A} : \tilde\rho _P (a) \neq  0 \} = \{P \in \operatorname{Prim} \mathcal{A} :\rho _{\tilde{P}} (a) \neq  0 \},
	\end{equation}
	where \(\tilde\rho _P \colon M(\mathcal{A}) \to M(\mathcal{A}/{P})\) is the extension to the multiplier algebra of the quotient map \(\rho _P\colon \mathcal{A} \to \mathcal{A}/{P}\), also form a base of the topology on \(\operatorname{Prim}(\mathcal{A})\). Here, the second equality follows from the relation \(\ker\tilde\rho_P=\tilde P\). To verify this relation, note that \(\mathcal{A}/P\) is an essential ideal in \(\tilde\rho_P(M(\mathcal{A}))\subset M(\mathcal{A}/P)\), while \(\{0\}\in\operatorname{Prim}(\mathcal{A}/P)\). Hence, \(\{0\}\in\operatorname{Prim}(\tilde\rho_P(M(\mathcal{A})))\), so \(\ker\tilde\rho_P\in\operatorname{Prim}(M(\mathcal{A}))\). Together with \((\ker\tilde\rho_P)\cap\mathcal{A}=\ker\rho_P=P\), this gives \(\ker\tilde\rho_P=\tilde P\).

	Let the point \((P,K)\) belong to the subspace \(\operatorname{Prim}^{\mathcal{A}}(\mathcal{A}/\mathcal{J}_0) \times \operatorname{Prim} (\mathcal{D}\rtimes_{\beta} G)\) and let \(U_b \times U_{d}\) be its open neighbourhood in \(\operatorname{Prim} \mathcal{A} \times \operatorname{Prim}(\mathcal{D}\rtimes _{\beta} G)\), where \(U_b\) and \(U_{d}\) are of the form \eqref{Eq:jacobson-top_base-system} for some \(b \in \mathcal{A}\) and \(d \in \mathcal{D}\rtimes _{\beta} G\). Put \(\eta = \Pi_1^{-1} \circ \phi \circ \rho_{\mathcal{J}_0} \colon \operatorname{Prim}^{\mathcal{A}}(\mathcal{A}/\mathcal{J}_0) \to \operatorname{Prim}\mathcal{C}\) and choose \(c \in \mathcal{C}\) such that \(\rho_{\eta(P)}(c) \neq 0\). Since \(\psi \circ (\phi \rtimes \operatorname{id}) \circ \rho_{\mathcal{J}_{0} \rtimes _\alpha G }\) is surjective, there exists \(\tilde{d} \in \mathcal{A}\rtimes _\alpha G\) such that 
	\[
	\psi \circ (\phi \rtimes \operatorname{id}) \circ \rho_{\mathcal{J}_{0} \rtimes _\alpha G }  (\tilde{d}) = c \otimes d \in \mathcal{C}\otimes (\mathcal{D}\rtimes _{\beta} G ).
	\]
	Let \(i_{\mathcal{A}}\) be the canonical embedding of \(\mathcal{A}\) into the multiplier algebra \(M(\mathcal{A}\rtimes _\alpha G)\), that is, for \(a \in \mathcal{A}\) and \(f\) in the dense subspace \(C_c(G,\mathcal{A})\) one has \((i_{\mathcal{A}}(a)f)(t) = a f(t)\) for all \(t \in G\) (see \cite[Proposition~2.34]{Will2007}).
	For any \(0 \neq z \in M(\mathcal{A}\rtimes _\alpha G)\) set \(a_z = i_{\mathcal{A}} (b) \cdot z \cdot  \tilde{d} \in M(\mathcal{A}\rtimes _\alpha G)\).
	We will show that the inclusion
	\[
	V_{a_z} \subset \Phi(U_{b} \times U_{d})
	\]
	holds, where \(V_{a_z}\) is of the form \eqref{Eq:jacobson-top_base-system-multiplier}. Let \(S\in V_{a_z}\) be a primitive ideal. From the inequalities
	\[
	0<\lVert \tilde\rho_{S} (a_z) \rVert \leq
	\lVert \tilde{\rho }_S(i_{\mathcal{A}} (b)) \rVert \cdot \lVert \tilde \rho_S (z) \rVert \cdot  \lVert \tilde{\rho
	}_S(\tilde{d}) \rVert
	\]
	it follows that \(\tilde{\rho }_S(i_{\mathcal{A}} (b)) \neq 0\) and \(\tilde{\rho }_S(\tilde{d}) \neq 0\). Since the dynamical system \((\mathcal{A},G,\alpha )\) is quasi-regular, there exists \(Q \in \mathcal{Q}_{\mathcal{A}}\) such that \(\operatorname{Res} S = \mathcal{I}_{Q}\). By \cite[Proposition~2.34]{Will2007}, from \(\tilde{\rho }_S(i_{\mathcal{A}} (b)) \neq 0\) we obtain \(\rho_{\mathcal{I}_Q}(b) \neq 0\). Hence, there exists a primitive ideal \(P' \in \kappa^{-1}(Q)\) with \(\rho_{P'}(b) \neq 0\), i.e. \(P' \in U_b\). Moreover, by surjectivity of \(\Phi\) there exists \(K' \in \operatorname{Prim} (\mathcal{D}\rtimes _\beta G)\) such that \(\Phi (P',K')=S\). If \(P' \in \operatorname{Prim}^{\mathcal{A}}\mathcal{J}_0\), then the map \(\Phi (P',\cdot)\) is constant, and \(S = \Phi (P',K) \in \Phi (U_b \times U_{d})\). If \(P' \in \operatorname{Prim}^{\mathcal{A}} \mathcal{A}/\mathcal{J}_0\), then from \(\tilde{\rho }_S(\tilde{d}) \neq 0\) and the equalities
	\[
	\ker \rho_{S} = \ker \rho_{\Phi (P',K')} = \ker \bigl(
	(\rho_{\eta(P')} \otimes  \rho_ {K'}) \circ \psi \circ (\phi \rtimes \operatorname{id}) \circ \rho_{\mathcal{J}_{0}\rtimes_\alpha G} \bigr)
	\]
	it follows that
	\[
	\bigl(\rho_{\eta(P')} \otimes \rho_ {K'}\bigr) \bigl(c \otimes d\bigr) = \rho_{\eta(P')}(c) \otimes \rho_ {K'}(d)\neq 0,
	\]
	i.e. \(\rho_{K'}(d) \neq 0\). Therefore, \(S = \Phi (P',K') \in \Phi (U_b \times U_{d})\) as well. Thus, for every \(z \in M(\mathcal{A}\rtimes _\alpha G)\) the inclusion \(V_{a_z} \subset\Phi (U_b \times U_{d})\) holds.
	
	Consider the open set \(V = \bigcup_{z \in M(\mathcal{A}\rtimes _\alpha G)} V_{a_z}\). Clearly, \(V \subset\Phi (U_b \times U_{d})\).
	We will show that \(\Phi (P,K) \in  V\). Suppose the contrary, i.e. assume that for every \(z \in M(\mathcal{A}\rtimes _\alpha G)\) the relation
	\(\tilde\rho_{\Phi (P,K)} (i_{\mathcal{A}} (b) \cdot z \cdot \tilde{d}) = 0\) holds.
	Since the ideal \( \left. \Phi (P,K) \right. \widetilde{} \,\) is primitive, it is also a prime ideal in \(M(\mathcal{A}\rtimes _\alpha G)\); hence it follows that either \(\tilde\rho_{\Phi (P,K)} (i_{\mathcal{A}} (b) ) = 0\) or
	\(\tilde\rho_{\Phi (P,K)} ( \tilde{d}) = 0\) (see \cite[II.5.4.5]{Blac2006}). In the first case, from \(\operatorname{Res} (\Phi (P,K)) = P\) (recall we assume \(\Phi (P,K) \in \operatorname{Prim}^{\mathcal{A} \rtimes _\alpha G} ((\mathcal{A}/\mathcal{J}_0)\rtimes _{\alpha^{\mathcal{J}_0}} G)\)) and \cite[Proposition~2.34]{Will2007} we obtain \(\rho _P(b) = 0\), i.e. \(P \notin U_b\). In the second case, from \(\bigl(\rho_{\eta(P)} \otimes \rho_ {K}\bigr) (c \otimes d) = \rho_{\eta(P)}(c) \otimes \rho_ {K}(d) = 0\) and \(\rho_{\eta(P)}(c) \neq 0\) it follows that \(\rho_ {K}(d) = 0\), i.e. \(K \notin  U_{d}\). Both possibilities contradict the original choice of \((P,K)\). Therefore, \(\Phi (P,K) \in  V \subset \Phi (U_b \times U_d)\), which proves that \(\Phi\) is open.

	The theorem is proved.
\end{proof}

\begin{remark}\label{rem:sim-rel-1}
	Consider the equivalence relation \(\sim\) on the space \(\operatorname{Prim} \mathcal{A} \times \operatorname{Prim} (\mathcal{D}\rtimes _{\beta} G)\), defined by \((P_1,K_1) \sim (P_2,K_2)\) if and only if \((P_1,K_1) = (P_2,K_2)\) when \(P_1, P_2 \in \operatorname{Prim} \mathcal{A}/\mathcal{J}_0\), or \(\kappa (P_1)= \kappa (P_2)\) when \(P_1, P_2 \in \operatorname{Prim} \mathcal{J}_0\). It is clear that if Theorem \ref{th:free-and-fixed-action} holds, then the map \(\tilde{\Phi } \left( [(P,K)] \right) = \Phi (P,K)\) is injective and therefore establishes a homeomorphism between the spaces \((\operatorname{Prim} \mathcal{A} \times \operatorname{Prim} (\mathcal{D}\rtimes _{\beta} G)) /_\sim \) and \(\operatorname{Prim} (\mathcal{A}\rtimes _\alpha G)\), with the equivalence relation \(\sim \) being open (cf. the paragraph after the proof of \cite[Theorem~5.3]{Will1981a}). Thus, starting from the topology of \(\operatorname{Prim} \mathcal{A} \times \operatorname{Prim} (\mathcal{D}\rtimes _{\beta} G)\), we obtain a description of the topology on the quotient space \((\operatorname{Prim} \mathcal{A} \times \operatorname{Prim} (\mathcal{D}\rtimes _{\beta} G)) /_\sim \) and on \(\operatorname{Prim} (\mathcal{A}\rtimes _\alpha G)\).

	In applications, the quotient space \((\operatorname{Prim} \mathcal{A} \times \operatorname{Prim} (\mathcal{D}\rtimes _{\beta} G)) /_\sim \) is sometimes more conveniently described in another way. The identification \(\sim\) can be performed in two steps. At the first step, the factorization occurs in the first component of the product \(\operatorname{Prim} \mathcal{A} \times \operatorname{Prim} (\mathcal{D}\rtimes _{\beta} G)\), replacing the spectrum \(\operatorname{Prim} \mathcal{A}\) with the quasi-orbit space \(\mathcal{Q}_{\mathcal{A}}\). At the second step, we consider the following equivalence relation \(\approx\): for \(Q_1, Q_2 \in \kappa (\operatorname{Prim} (\mathcal{A}/\mathcal{J}_0)) \) we have \((Q_1, K_1) \approx (Q_2,K_2)\) if and only if \(Q_1 = Q_2\) and \(K_1=K_2\), while for \(Q_1, Q_2 \in \kappa (\operatorname{Prim} \mathcal{J}_0)\) the relation \((Q_1,K_1) \approx (Q_2,K_2)\) holds if and only if \(Q_1=Q_2\). Then, the spaces \((\operatorname{Prim} \mathcal{A} \times \operatorname{Prim} (\mathcal{D}\rtimes _{\beta} G)) /_\sim \) and \((\mathcal{Q}_{\mathcal{A}} \times \operatorname{Prim} (\mathcal{D}\rtimes _{\beta} G)) /_\approx\) coincide.
\end{remark}

\begin{remark}
	The amenability condition on the group \(G\) in Theorem \ref{th:free-and-fixed-action} can be replaced by weaker assumptions. It is sufficient to require that the so-called reduced crossed product \(\mathcal{A} \rtimes_{\alpha,r} G \) (see \cite[Definition~7.7]{Will2007}) coincide with the full crossed product \(\mathcal{A} \rtimes_\alpha G\) (more precisely, we require the existence of an isomorphism that agrees on the dense subset \(C_c(G,\mathcal{A})\)), and that the group \(G\) be exact (see \cite[the definition after (7.11)]{Will2007}).

	Indeed, in the proof of Theorem \ref{th:free-and-fixed-action} the amenability of the group \(G\) was used only in reference to Proposition \ref{prop:prim_is_fibered:opennes} and Lemma \ref{lm:Phi-on-prim-J0}. In the proofs of these statements, this property of the group was used only to deduce the equalities \(\operatorname{Ind} \mathcal{I} = \mathcal{I} \rtimes_\gamma G\) for various dynamical systems \((\mathcal{B},G,\gamma)\) and invariant ideals \(\mathcal{I}\) in the algebra \(\mathcal{B}\).
	The dynamical systems \((\mathcal{B},G,\gamma)\) arising there are obtained from the initial system \((\mathcal{A},G,\alpha)\) by passing to systems associated with invariant ideals or quotient algebras by invariant ideals. Therefore, to ensure the validity of Theorem \ref{th:free-and-fixed-action}, it is sufficient that the equalities \(\operatorname{Ind} \mathcal{I} = \mathcal{I} \rtimes_\gamma G\) hold for all such systems.
	The equality \(\operatorname{Ind} \mathcal{I} = \mathcal{I} \rtimes_\gamma G\) holds if and only if \((\mathcal{B} /\mathcal{I}) \rtimes_{\gamma} G \cong (\mathcal{B} /\mathcal{I}) \rtimes_{\gamma,r} G \) (\cite[Remark~5.24]{Will2007}). Hence, it suffices to show that for any dynamical system \((\mathcal{B},G,\gamma)\), the isomorphism \(\mathcal{B} \rtimes_\gamma G \cong \mathcal{B} \rtimes_{\gamma,r} G \) implies
	\begin{equation}\label{eq:reduced=full:ideal-quotient}
		\mathcal{I}\rtimes_\gamma G \cong  \mathcal{I} \rtimes_{\gamma,r} G \quad \text{and} \quad (\mathcal{B} /\mathcal{I}) \rtimes_{\gamma^\mathcal{I}} G  \cong (\mathcal{B}/\mathcal{I}) \rtimes_{\gamma^{\mathcal{I}},r} G.
	\end{equation}
	From the inclusions \(C_c(G,\mathcal{I}) \subset C_c(G,\mathcal{B})\) and \(\mathcal{I} \rtimes_{\gamma,r} G \subset \mathcal{B} \rtimes_{\gamma,r} G\) we obtain the first relation in \eqref{eq:reduced=full:ideal-quotient}. From this and the exactness of \(G\) we have
	\[
		(\mathcal{B}/\mathcal{I}) \rtimes_{\gamma^\mathcal{I},r} G  = (\mathcal{B} \rtimes _{\gamma,r} G)/(\mathcal{I} \rtimes _{\gamma,r} G) \cong (\mathcal{B} \rtimes _\gamma G)/(\mathcal{I} \rtimes _\gamma G) = (\mathcal{B}/\mathcal{I}) \rtimes _{\gamma^{\mathcal{I}} } G,
	\]
	as required.
\end{remark}

\appendix
\section{Corrections to the published version}\label{app:corrections}

This is a preprint of the Work accepted for publication in Mathematical Notes, \copyright~2026, \href{https://pleiades.online/}{Pleiades Publishing, Ltd.} Published version: \href{https://doi.org/10.1134/S0001434626602030}{Math. Notes \textbf{119}~(4), 559--574 (2026)}. Russian original: \href{https://doi.org/10.4213/mzm14787}{Mat. Zametki \textbf{119}~(4), 471--491 (2026)}. Translated from the Russian by the author.

The numbering of all statements and formulas is preserved. In this version, we fix minor misprints and add a self-contained proof of Lemma~\ref{lm:locally-closed-R-Q} (previously only a reference was given). Finally, we correct two inaccuracies in the proof of Theorem~\ref{th:free-and-fixed-action}, leaving its statement unchanged.

Specifically, the proof of Theorem~\ref{th:free-and-fixed-action} is modified as follows. First, since the isomorphism \(M(\mathcal{A}/P)\cong M(\mathcal{A})/\tilde P\) does not hold in general, the second equality in~\eqref{Eq:jacobson-top_base-system-multiplier} is justified directly by the relation \(\ker\tilde\rho_P=\tilde P\). Second, since the map \(\tilde\rho_{\mathcal{J}_0\rtimes_\alpha G}\) need not be surjective, instead of lifting \(1\otimes d\) to \(M(\mathcal{A}\rtimes_\alpha G)\) we take \(\tilde d\in\mathcal{A}\rtimes_\alpha G\) satisfying \(\psi\circ(\phi\rtimes\operatorname{id})\circ\rho_{\mathcal{J}_0\rtimes_\alpha G}(\tilde d)=c\otimes d\in\mathcal{C}\otimes(\mathcal{D}\rtimes_\beta G)\) for an element \(c\in\mathcal{C}\) with \(\rho_{\Pi_1^{-1}(\phi(\rho_{\mathcal{J}_0}(P)))}(c)\neq 0\).

\end{document}